\documentclass{article}
	
\usepackage{amsthm,amsfonts,amsmath,amscd,amssymb,latexsym,epsfig}
\usepackage{mathrsfs}
\usepackage{mathdots}
\usepackage[all]{xy}
\usepackage{hyperref}

\usepackage {dsfont}
\usepackage{enumitem}

\usepackage[T1]{fontenc}
\usepackage[utf8]{inputenc}
	
\usepackage[a4paper,margin=1.5in]{geometry}

	\newtheorem{TheoremABC}{Theorem}

	\newtheorem {Theorem}{Theorem}[section]
	\newtheorem {Cor}[Theorem]{Corollary}
	\newtheorem {Lemma}[Theorem]{Lemma}
	\newtheorem {Proposition}[Theorem]{Proposition}
	\newtheorem {Definition}[Theorem]{Definition}

	\theoremstyle{definition}
	\newtheorem {Remark}[Theorem]{Remark}

\begin {document}

\title{Convergence of the Yang--Mills--Higgs flow on gauged holomorphic maps and applications\thanks{The author was partially supported by the Swiss National
Science Foundation (grant number  200021-156000).}}

\author{Samuel Trautwein\thanks{samuel.trautwein@math.ethz.ch, ETH Zürich, Rämistrasse 101, 8902 Zürich (Switzerland)}}

\maketitle

\abstract{The symplectic vortex equations admit a variational description as global minimum of the Yang--Mills--Higgs functional. We study its negative gradient flow on holomorphic pairs $(A,u)$ where $A$ is a connection on a principal $G$-bundle $P$ over a closed Riemann surface $\Sigma$ and $u:  P \rightarrow X$ is an equivariant map into a Kähler Hamiltonian $G$-manifold. The connection $A$ induces a holomorphic structure on the Kähler fibration $P\times_G X$ and we require that $u$ descends to a holomorphic section of this fibration.

We prove a {\L}ojasiewicz type gradient inequality and show uniform convergence of the negative gradient flow in the $W^{1,2}\times W^{2,2}$-topology when $X$ is equivariantly convex at infinity with proper moment map, $X$ is holomorphically aspherical and its Kähler metric is analytic. 

As applications we establish several results inspired by finite dimensional GIT: First, we prove a certain uniqueness property for the critical points of the Yang--Mills--Higgs functional which is the analogue of the Ness uniqueness theorem. Second, we extend Mundet's Kobayashi--Hitchin correspondence to the polystable and semistable case. The arguments for the polystable case lead to a new proof in the stable case. Third, in proving the semistable correspondence, we establish the moment-weight inequality for the vortex equation and prove the analogue of the Kempf existence and uniqueness theorem.}

\newpage

\setcounter{tocdepth}{2}
\tableofcontents

\newpage

\section{Introduction}
The symplectic vortex equations \cite{Salamon:2000, Salamon:2003, Mundet:2000} are an equivariant version of the $J$-holomorphic curves equation in symplectic geometry. These equations also generalize the Yang--Mills equations \cite{AtBott:YangMillsEq}, the notion of Bradlow pairs \cite{Bradlow:1991} and are closely related to Hitchin's selfduality equations \cite{Hitchin:1987} and Higgs-bundles. 

Solutions of the symplectic vortex equation admit a variational characterization as minima of the Yang--Mills--Higgs functional. The main analytic result of this paper is a {\L}ojasiewicz type gradient inequality for this functional and the convergence of the associated negative gradient flow under suitable assumptions (Theorem \ref{ThmConv0}). As applications we obtain several new results motivated by geometric invariant theory. These are a uniqueness result for critical points of the Yang--Mills--Higgs functional (Theorem \ref{ThmB}), an analogue of the Kempf--Ness Theorem in the present setting (Theorem \ref{ThmC}), an extension of Mundet's Kobayashi-Hitchin correspondence (Theorem \ref{ThmD}) and a sharp moment weight inequality together with the existence and uniqueness theorem for the dominant weight (Theorem \ref{ThmE}). Our proofs are inspired by the work of Calabi, Chen, Donaldson, Sun \cite{CalabiChen:2002, ChenSun:2010, Chen:2008, Chen:2009, Donaldson:2005} on extremal Kähler metrics; see \cite{RobSaGeo} for a finite dimensional discussion.

\textbf{Setup.} Throughout this paper we assume the following. $G$ is a compact (real) Lie group with Lie algebra $\mathfrak{g}$ together with a fixed choice of an invariant inner product on $\mathfrak{g}$, $\Sigma$ is a closed Riemann surface with fixed volume form $dvol_{\Sigma}$ and induced Riemannian metric, $P \rightarrow \Sigma$ is a principal $G$ bundle and $(X,J, \omega)$ is a Kähler manifold equipped with a Hamiltonian $G$ action induced by an equivariant moment map $\mu: X \rightarrow \mathfrak{g}$.

\subsection{Geometric invariant theory for the vortex equation}
Atiyah--Bott \cite{AtBott:YangMillsEq} observed that the curvature $F_A \in \Omega^2(\Sigma, \text{ad}(P))$ defines a moment map for the action of the gauge group $\mathcal{G}(P)$ on the space of connections $\mathcal{A}(P)$. The vortex equations are obtained as an extension of this picture. Consider the associated Kähler fibration 
									$$P(X) := P\times_G X := (P\times X)/G.$$ 
and denote by $\mathcal{S}(P,X)$ its space of sections. The symplectic vortex equations on pairs $(A,u) \in \mathcal{A}(P)\times\mathcal{S}(P,X)$ are given by
			\begin{align} \label{vortexeq1} \bar{\partial}_A u = 0, \qquad *F_A + \mu(u) = 0. \end{align}
The connection $A \in \mathcal{A}(P)$ induces a holomorphic structure on the total space of the Kähler fibration $P(X)$ and the equation $\bar{\partial}_A u = 0$ requires $u$ to be a holomorphic section. The subspace
					$$\mathcal{H}(P,X) := \{ (A,u) \in \mathcal{A}(P)\times \mathcal{S}(P,X)\,|\, \bar{\partial}_A u = 0\}$$
is formally a Kähler submanifold of $\mathcal{A}(P)\times\mathcal{S}(P,X)$. It is well known that 
			\begin{align} \label{introeq0} \Phi: \mathcal{A}(P)\times\mathcal{S}(P,X) \rightarrow \Omega^0(\Sigma, \text{ad}(P)), \qquad \Phi(A,u) := *F_A + \mu(u) \end{align}
provides a moment map for the $\mathcal{G}(P)$-action on $\mathcal{H}(P,X)$ (see Lemma \ref{LemmaMoment}) and solutions of (\ref{vortexeq1}) give rise to the symplectic moduli space
				$$\mathcal{M}_{symp}(P,X) := \left\{ (A,u) \in \mathcal{H}(P,X)\,|\, *F_A + \mu(u) = 0\right\}/\mathcal{G}(P).$$

This moduli space admits an alternative description as complex GIT quotient of $\mathcal{H}(P,X)$. For this let $G^c$ be the complexification of $G$, let $P^c := P \times_G G^c$ be the complexification of $P$ and define the complexified gauge group as $\mathcal{G}^c(P) := \mathcal{G}(P^c)$. There exists a one to one correspondence between smooth connections on $P$ and holomorphic structures on $P^c$ (see \cite{Singer:1959}). This yields a natural action of $\mathcal{G}^c(P)$ on $\mathcal{A}(P)$ which extends the gauge action. Assume that the $G$-action on $(X,J,\omega)$ extends to a holomorphic $G^c$-action on $(X,J)$ such that $\mathcal{G}^c(P)$ acts naturally on $\mathcal{S}(P,X)$. 


\begin{Definition} \label{DefStability}
Let $(A,u)\in \mathcal{H}(P,X)$ and denote by $\overline{\mathcal{G}^c(A,u)}$ the $W^{1,2}\times W^{2,2}$-closure of its complexified orbit\footnote{
Here it suffices to consider the closure within the space $\mathcal{H}(P,X)$ of smooth holomorphic pairs. In the main part of the paper we will consider pairs $(A,u)$ of Sobolev class $W^{1,2}\times W^{2,2}$ and gauge transformations of Sobolev class $W^{2,2}$. This does not affect the overall picture since (a) every complex orbit contains a dense set of smooth representatives and (b) every $W^{1,2}\times W^{2,2}$ solution to the vortex equation is gauge equivalent to a smooth solution. See Section \ref{Section.SobolevReg} and Lemma \ref{Lemma.SobolevReg}.
}. Denote by $\Phi(A,u) := *F_A + \mu(u)$ the moment map (\ref{introeq0}).
		\begin{enumerate}
			\item	$(A,u)$ is called \textbf{stable}, if $\Phi^{-1}(0)\cap \mathcal{G}^c(A,u) \neq \emptyset$ and the isotropy subgroup $\mathcal{G}_{(A,u)} := \{ k \in \mathcal{G}(P)\,|\, k(A,u) = (A,u)\}$ is discrete.
			\item	$(A,u)$ is called \textbf{polystable}, if $\Phi^{-1}(0)\cap \mathcal{G}^c(A,u) \neq \emptyset$.
			\item	$(A,u)$ is called \textbf{semistable}, if $\Phi^{-1}(0)\cap \overline{\mathcal{G}^c(A,u)} \neq \emptyset$ .
			\item	$(A,u)$ is called \textbf{unstable}, if $\Phi^{-1}(0)\cap \overline{\mathcal{G}^c(A,u)} = \emptyset$.
		\end{enumerate}
Denote by $\mathcal{H}^s \subset \mathcal{H}^{ps} \subset \mathcal{H}^{ss}$ and $\mathcal{H}^{us}$ the corresponding $\mathcal{G}^c(P)$-invariant subspaces.
\end{Definition}

The GIT quotient of $\mathcal{H}(P,X)$ by $\mathcal{G}^c(P)$ is defined as the quotient space
				$$\mathcal{M}_{GIT}(P,X) := \mathcal{H}^{ss}(P,X)/\!/\mathcal{G}^c(P) := (\mathcal{H}^{ss}(P,X)/\mathcal{G}^c(P))/\sim$$
under the orbit closure relation $\mathcal{G}^c (A,u) \sim \mathcal{G}^c (B,v)$ if and only if $\overline{\mathcal{G}^c (A,u)} \cap \overline{\mathcal{G}^c (B,v)} \cap \mathcal{H}^{ss}(P,X) \neq \emptyset$. It follows from our main results that each equivalence class in this quotient contains a unique $\mathcal{G}(P)$-orbit of solutions to the symplectic vortex equations and $\mathcal{M}_{GIT}(P,X) \cong \mathcal{M}_{symp}(P,X)$ (see Corollary \ref{CorMain}).

\subsection{The main theorem}

The moment map squared functional plays a crucial role in the differential geometric version of GIT. It is defined by
			\begin{align} \label{DefF} \mathcal{F}: \mathcal{H}(P,X) \rightarrow \mathbb{R}, \qquad \mathcal{F}(A,u) := \frac{1}{2}\int_{\Sigma} ||*F_A + \mu(u)||^2 \, dvol_{\Sigma} \end{align}
and closely related to the Yang--Mills--Higgs functional
			\begin{align} \label{DefYMH} \mathcal{YMH}(A,u) := \frac{1}{2} \int_{\Sigma} ||F_A||^2 + ||d_A u||^2 + ||\mu(u)||^2 \, dvol_{\Sigma} \end{align}
by the energy identity in Proposition \ref{PropEngergy}. In particular, for $(A,u) \in \mathcal{H}(P,X)$ it holds $\nabla \mathcal{YMH}(A,u) = \nabla \mathcal{F}(A,u)$, albeit the gradients look quite different at first glance. The negative gradient flow on $\mathcal{H}(P,X)$ has the following form
			\begin{equation} \label{floweq0}  
				\begin{gathered} 
							A(0) = A_0, \qquad u(0) = u_0, \qquad \bar{\partial}_A(u) = 0\\
										\partial_t A = -*d_A( *F_A + \mu(u)), \qquad 		\partial_t u =   J L_u(*F_A + \mu(u)) 
				\end{gathered}
			\end{equation}							
Our main result says that solutions exist for all time and converge under the following hypothesis:
			
		\begin{itemize}	
			\item[\textbf{(A)}]  The Kähler metric on $X$ and the moment map $\mu: X \rightarrow \mathfrak{g}$ are both analytic.
			\item[\textbf{(B)}]  $X$ is holomorphically aspherical.
			\item[\textbf{(C)}]  $\mu$ is proper and $X$ is equivariantly convex at infinity, i.e. there exists a proper $G$-invariant function $f: X \rightarrow [0,\infty)$ and $c_0 > 0$ such that
				\begin{align} \label{EqC} f(x) \geq c_0 \quad \Longrightarrow\quad  \begin{array}{c} \langle \nabla_{\xi} \nabla f(x), v \rangle + \langle \nabla_{Jv} \nabla f(x), J v \rangle \geq 0 \\ df(x)JL_x\mu(x) \geq 0 \end{array} \end{align}
			for every $x \in X$	and $v \in T_x X$.
		\end{itemize}

\begin{TheoremABC}[\textbf{Convergence}] \label{ThmConv0}
Assume \textbf{(C)} and let $(A_0,u_0) \in \mathcal{H}(P,X)$ be given. Then there exists a unique solution 
				$$(A,u): [0,\infty) \rightarrow \mathcal{H}(P,X)$$ 
of (\ref{floweq0}) which exists for all times $t \geq 0$. If in addition \textbf{(A)}, \textbf{(B)} are satisfied, then there exists a critical point $(A_{\infty}, u_{\infty}) \in \mathcal{A}^{1,2}(P)\times \mathcal{S}^{2,2}(P,X)$ of Sobolev class $W^{1,2}\times W^{2,2}$ and $T, C, \epsilon > 0$ such that for all $t > T$ the pointwise distance between $u(t)$ and $u_{\infty}$ is smaller then the injectivity radius of $X$ along $u_{\infty}(P)$ and
					$$||A(t) - A_{\infty}||_{W^{1,2}} + ||\exp_{u_{\infty}}^{-1}  u(t)||_{W^{2,2}} \leq C t^{-\epsilon}.$$
\end{TheoremABC}

\begin{proof}
In Theorem \ref{ThmVen} the existence is proven together with certain continuity and regularity assertions on the flow. The convergence part is proven in Theorem \ref{ThmConv}.
\end{proof}

\begin{Remark}[\textbf{Regularity of the Limit.}]
Starting at a smooth initial condition $(A_0,u_0) \in \mathcal{H}(P,X)$, the solution $(A(t),u(t))$ of (\ref{floweq0}) remains smooth for all times $t > 0$. However, it is an open question if the limit $(A_{\infty},u_{\infty})$ is smooth.
\end{Remark}

Lin \cite{Lin:2012} and Venugopalan \cite{Venugopalan:2016} discussed the flow (\ref{floweq0}) independently and they proved under certain hypotheses that solutions exist for all times. Lin \cite{Lin:2012} considered in fact a generalization of (\ref{floweq0}), where $\Sigma$ is replaced by a compact Kähler manifold, and showed that smooth solutions exist for all times when $X$ is compact. His proof follows ideas of Donaldson \cite{Donaldson:ASD4} and he translates (\ref{floweq0}) into a heat flow on the space of complex gauge connections. Venugopalan \cite{Venugopalan:2016} extended the arguments given by R\r{a}de \cite{Rade1992} for the Yang--Mills flow and proved short time existence together with an uniform lower bound of the existence interval. For this argument she needed to assume that the flow remains in a compact region of $X$. We verify in Lemma \ref{LemmaEC} that this property follows from \textbf{(C)} and the maximum principle.

The main ingredient in our proof of the convergence of solutions to (\ref{floweq0}) is a {\L}ojasiewicz gradient inequality for the Yang--Mills--Higgs functional (Theorem \ref{ThmLoj1}). This approach was introduced by Simon \cite{Simon:1983} and in its implementation we follow the arguments given by R\r{a}de \cite{Rade1992} for the Yang--Mills flow.

\begin{Remark}[\textbf{On assumption (A)}]
The proof depends on a suitable version of the {\L}ojasiewicz gradient inequality and requires an analytic setup. In the finite dimensional case, it follows from the Marle and Guillemin-Sternberg normal form that the moment map squared functional is locally analytic (see Lerman \cite{Lerman:2005}). If an analogous result is valid in our infinite dimensional setting, one might hope to remove this assumption.
\end{Remark}

\begin{Remark}[\textbf{On assumption (B)}] \label{RemarkB}
\begin{enumerate}
	\item Holomorphically aspherical means that every holomorphic map $\mathbb{C}P^1 \rightarrow X$ is constant.
	\item This assumption prevents bubbling of holomorphic spheres within the fiber and is needed to establish sequential compactness along the flow lines (see Proposition \ref{PropCpct}).
	
	\item $X$ is necessarily noncompact under this assumption. Suppose otherwise that $X$ is compact and there exists $\xi \in \mathfrak{g}\backslash\{0\}$ and $x_0 \in X$ such that $\exp(\xi) = \mathds{1}$ and the infinitesimal action $L_{x_0}\xi \neq 0$ is nontrivial. Let $x: \mathbb{R} \rightarrow X$ be the solution of $\dot{x} = -J L_x \xi = -\nabla H_{\xi}(x)$ with $H_{\xi}:= \langle \mu, \xi \rangle$ starting at $x(0) = x_0$. Since $H_{\xi}$ is a Morse--Bott function, $x(t)$ converges exponentially to critical points $x^{\pm}$ as $t \rightarrow \pm \infty$ satisfying $L_{x^{\pm}} \xi = 0$. Using the $S^1$ action obtained from integrating the infinitesimal action of $\xi$, one can rotate this flow line within $X$ and construct a nontrivial holomorphic sphere.
	
	\item When $X$ has nonpositive curvature, the distance function is plurisubharmonic and every holomorphic sphere $\mathbb{C}P^1 \rightarrow X$ is constant.
\end{enumerate}
\end{Remark}

\begin{Remark}
Important examples in which our assumptions are satisfied arise when $X$ is a complex vector space (see \cite{Banfield:1999, BrDaskGPWent:1995}).
\end{Remark}

\begin{Remark} [\textbf{Higgs bundles}]
Let $X = \mathfrak{g}^c$ and consider the adjoint action of $G$ on $\mathfrak{g}^c$. This action is Hamiltonian with moment map $\mu(\zeta) = \frac{i}{2}[\zeta, \zeta^*]$ where $\zeta^* := -\text{Re}(\zeta) + \textbf{i}\text{Im}(\zeta)$. Then $P(X) = \text{ad}(P^c)$ is a holomorphic vector bundle and our assumptions are satisfied. Higgs bundles are obtained as a slight variant of this setup where one considers holomorphic sections of the twisted bundle $P(X)\otimes K := \Omega^{1,0}(\Sigma, \text{ad}(P^c))$. While this is not covered by our general discussion, the proof generalizes ad verbatim to this case.
\end{Remark}

\subsection{Consequences of the main theorem}

The infinitesimal action of $\xi \in \Omega^0(\Sigma, \text{ad}(P))$ on $\mathcal{A}(P)\times \mathcal{S}(X,P)$ is given by
	\begin{align*}
			\mathcal{L}_{(A,u)} \xi := \left.\frac{d}{dt}\right|_{t=0} \exp(t\xi) (A,u) = \left( -d_A \xi, L_u \xi \right)
	\end{align*}
where $L_x : \mathfrak{g} \rightarrow T_x X$ denotes the infinitesimal action of $\mathfrak{g}$ on $X$. Denote by $\mathcal{L}^c_{(A,u)}$ the infinitesimal action of $\mathcal{G}^c(P)$ which agrees with the complexification of $\mathcal{L}_{(A,u)}$. Then 
				$$\nabla \mathcal{F} (A,u) = \mathcal{L}_{(A,u)}^c \textbf{i} (*F_A + \mu(u)) $$
implies that solutions of (\ref{floweq0}) remain in a single complexified orbit. The following result is the analogue of the Ness uniqueness theorem in finite dimensional GIT.
 
\begin{TheoremABC}[\textbf{Uniqueness of critical points}] \label{ThmB}
Assume \textbf{(A)}, \textbf{(B)} and \textbf{(C)}. 
	\begin{enumerate}
			\item Let $(A_0,u_0) \in \mathcal{H}(P,X)$ and let $(A_{\infty},u_{\infty})$ be the limit of the gradient flow (\ref{floweq0}) starting at $(A_0,u_0)$. Then 
					$$ ||*F_{A_{\infty}} + \mu(u_{\infty})||_{L^2} = \inf_{g \in \mathcal{G}^c(P)} ||*F_{gA_{0}} + \mu(g u_{0})||_{L^2} =:m.$$
			\item Suppose $(B_0,v_0), (B_1,v_1) \in \overline{\mathcal{G}^c (A_0,u_0)}$ (the $W^{1,2}\times W^{2,2}$-closure) and 
								$$||*F_{B_0} + \mu(v_0)||_{L^2} = m = ||*F_{B_1} + \mu(v_1)||_{L^2}.$$
			Then there exists $k \in \mathcal{G}(P)$ such that $(B_1,v_1) = k(B_0,v_0)$.
	\end{enumerate}						
\end{TheoremABC}

\begin{proof}
This is reformulated and proven in Theorem \ref{ThmNU}. The equivalence of both formulations follows from Proposition \ref{PropEngergy}.
\end{proof}

\begin{Cor} \label{CorMain}
Assume \textbf{(A)}, \textbf{(B)} and \textbf{(C)}. Every semistable orbit contains a unique polystable orbit in its $W^{1,2}\times W^{2,2}$-closure and every polystable orbit contains a unique $\mathcal{G}(P)$-orbit of solutions to the symplectic vortex equations.
\end{Cor}

The corollary shows $\mathcal{M}_{symp}(P,X) \cong \mathcal{M}_{GIT}(P,X)$. More explicitly, this isomorphism is obtained by the map which sends $(A_0,u_0) \in \mathcal{H}(P,X)$ to its limit $(A_{\infty}, u_{\infty})$ under (\ref{floweq0}). Theorem \ref{ThmLimitStability} gives a complete characterization for the different stability conditions in Definition \ref{DefStability} in terms of the limit $(A_{\infty}, u_{\infty})$.\\

Next, we need to recall the general construction behind the Kempf--Ness theorem. Given $(A,u) \in \mathcal{H}(P,X)$ there exists a $\mathcal{G}(P)$-invariant functional
				$$\Psi_{(A,u)} : \mathcal{G}^c(P) \rightarrow \mathbb{R}$$
whose gradient flow intertwines with (\ref{floweq0}) under the map $g \mapsto g^{-1}(A,u)$. The Kempf--Ness theorem characterizes the stability conditions of $(A,u)$ in Definition \ref{DefStability} in terms of the global properties of $\Psi_{(A,u)}$. The stable case is the main step in Mundet's proof of the Kobayashi--Hitchin correspondence \cite{Mundet:2000} and relates the stability of $(A,u)$ to a certain properness of $\Psi_{(A,u)}$. The remaining cases are the content of the next theorem, whose proof is a relatively easy consequence of Theorem \ref{ThmConv0} and Theorem \ref{ThmB}.

\begin{TheoremABC}[\textbf{Kempf--Ness Theorem}] \label{ThmC}
Assume \textbf{(A)}, \textbf{(B)}, \textbf{(C)} and let $(A,u) \in \mathcal{H}(P,X)$.
		\begin{enumerate}
				\item $(A,u)$ is polystable if and only if $\Psi_{(A,u)}$ has a critical point.
				\item $(A,u)$ is semistable if and only if $\Psi_{(A,u)}$ is bounded below.
				\item $(A,u)$ is unstable if and only if $\Psi_{(A,u)}$ is unbounded below.								
		\end{enumerate}
\end{TheoremABC}

\begin{proof} This is established in Theorem \ref{ThmKN}. \end{proof}

The weights for the $\mathcal{G}^c(P)$-action are defined as the asymptotic slopes of $\Psi_{(A,u)}$ along geodesics rays in $\mathcal{G}^c(P)/\mathcal{G}(P)$. For $(A,u) \in \mathcal{H}(P,X)$ and $\xi \in \Omega^0(\Sigma, \text{ad}(P))$ one has the explicit description
				$$w((A,u),\xi) := \lim_{t\rightarrow \infty} \left\langle *F_{e^{\textbf{i}t\xi}A} + \mu(e^{\textbf{i}t\xi} u), \,\xi \right\rangle_{L^2} \in \mathbb{R}\cup \{\infty\}.$$
Mundet's Kobayashi--Hitchin correspondence asserts that $(A,u)$ is stable if and only if $w((A,u), \xi) > 0$ for all $\xi \neq 0$. We extend this correspondence to the polystable and semistable case under the technical assumption on the pair $(A,u) \in \mathcal{H}(P,X)$.
		\begin{enumerate}[leftmargin=12mm]
				\item[\textbf{(H)}] For all $\xi \in \Omega^0(\Sigma, \text{ad}(P))$ with $w((A,u),\xi) \leq 0$ it holds $\sup_{t>0} ||\mu(e^{\textbf{i}t\xi} u)||_{L^2} < \infty.$
		\end{enumerate}					

\begin{Remark}[\textbf{On assumption (H)}] \label{Rmk.H}
\begin{enumerate}
	\item \textbf{(H)} is trivially satisfied for stable pairs $(A,u)$ and, by Proposition \ref{PropPS}, it is always satisfied for polystable pairs.
	\item By Proposition \ref{PropFW}, $w((A,u),\xi) < \infty$ implies that $A_+ := \lim_{t\rightarrow\infty} e^{\textbf{i}t\xi}A$ exists in $C^{\infty}$.

	\item	Proposition \ref{PropFW} provides a strong tool to verify \textbf{(H)}. When $X$ is a unitary vector space with linear $G \subset U(n)$ action, one can show that
					$$w((A,u),\xi) < \infty \qquad \Longrightarrow \qquad \lim_{t\rightarrow \infty} e^{\textbf{i}t\xi}(A,u) =: (A_+,u_+) $$
	where the limit exists in $C^{\infty}$ and \textbf{(H)} is satisfied in this case. Similarly, using Proposition \ref{PropFW}, one verifies \textbf{(H)} for Higgs bundles.

	\item \textbf{(H)} admits the following geometric description: For $(A,u) \in \mathcal{H}(P,X)$ denote by $\Psi_{(A,u)}: \mathcal{G}^c(P) \rightarrow \mathbb{R}$ its Kempf--Ness functional. 
					\begin{enumerate}[leftmargin=12mm]
				\item[\textbf{(H')}] For all $\xi \in \Omega^0(\Sigma, \text{ad}(P))$ the following holds: If  $\sup_{t > 0} \Psi_{(A,u)}\left(e^{-\textbf{i}t\xi}\right) < \infty$, then $\sup_{t>0} \left|\left|\nabla \Psi_{(A,u)}\left( e^{-\textbf{i}t\xi} \right) \right|\right|_{L^2} < \infty$.
		\end{enumerate}					
	Unraveling the definitions shows $\textbf{(H)} \Leftrightarrow \textbf{(H')}$. This property is reasonable to expect, since $\Psi_{(A,u)}$ is convex along geodesics. However, one can construct examples which show that convexity of $\Psi_{(A,u)}$ alone does not guarantee \textbf{(H')}.
	\item Unfortunately, we know little about the validity of \textbf{(H)} in general: We could neither prove that it is always satisfied, nor construct an explicit counterexample. This question is already meaningful (and open) in the finite dimensional case where $\Sigma = \{pt\}$.
\end{enumerate}
\end{Remark}
		
Consider the following properties for a pair $(A,u) \in \mathcal{H}(P,X)$:
 \begin{enumerate}[leftmargin=12mm]
				\item[\textbf{(SS)}] For all $\xi \in \Omega^0(\Sigma,\text{ad}(P))$ it holds $w((A,u),\xi) \geq 0$.
		
				\item[\textbf{(PS)}] For all $\xi \in \Omega^0(\Sigma,\text{ad}(P))$ with $\exp(\xi) = \mathds{1}$ and $w((A,u),\xi) = 0$ the limit 
								$\lim_{t\rightarrow \infty} e^{\textbf{i}t\xi}(A,u) \in (\mathcal{G}^c)^{2,2}(A,u)$
							exists in $W^{1,2}\times W^{2,2}$ and remains in the Sobolev completion of the complex group orbit.
		\end{enumerate}

\begin{TheoremABC}[\textbf{Polystable and semistable correspondence}] \label{ThmD}
Assume \textbf{(A)}, \textbf{(B)}, \textbf{(C)}, \textbf{(H)} and let $(A,u) \in \mathcal{H}(P,X)$.
		\begin{enumerate}
				\item $(A,u)$ is polystable if and only if it satisfies \textbf{(SS)} and \textbf{(PS)}.
					
				\item $(A,u)$ is semistable if and only if it satisfies \textbf{(SS)}.
		\end{enumerate}
\end{TheoremABC}

\begin{proof} This is established in Theorem \ref{ThmPS} and Theorem \ref{ThmSS}. \end{proof}

The polystable case has been established for twisted Higgs-bundles over Riemann surface by Garc\'{\i}a-Prada, Gothen and Mundet \cite{Garcia:2009} by different methods. They construct a Jordan-Hölder reduction and then deduce the polystable case from the stable case. For our proof the opposite is true and the stable case can be recovered as a special case of the polystable case. The proof is based on arguments of Chen--Sun \cite{ChenSun:2010}.

The semistable correspondence follows from a sharp version of the moment weight inequality stated next.

\begin{TheoremABC}[\textbf{Sharp moment-weight inequality}] \label{ThmE}
Assume \textbf{(H)}. For all $(A,u) \in \mathcal{H}(P,X)$ and $\xi \in \Omega^0(\Sigma, \text{ad}(P))\backslash\{0\}$ it holds
		\begin{align} \label{ThmEeq}  -\frac{w((A,u), \xi)}{||\xi||_{L^2}} \leq \inf_{g \in \mathcal{G}^c(P)} ||*F_{gA} + \mu(gu)||_{L^2}. \end{align}
If in addition \textbf{(A)}, \textbf{(B)}, \textbf{(C)} are satisfied and the right hand side is positive, then there exists a unique $\xi_0 \in \Omega^0(\Sigma, \text{ad}(P))$ with $||\xi_0||_{L^2} = 1$ which yields equality.
\end{TheoremABC}

\begin{proof} This is established in Theorem \ref{ThmSMWI}. \end{proof}

For finite dimensional projective spaces the estimate (\ref{ThmEeq}) is due to Mumford \cite{GIT:book} and Ness \cite[Lemma 3.1]{Ness:1984}, and the existence of a dominant weight is due to Kempf \cite{Kempf:1978}. Around the same time Atiyah--Bott \cite{AtBott:YangMillsEq} established this result for the Yang--Mills equations over Riemann surfaces. Its generalization to the hermitian Yang--Mills equations over higher dimensional base manifolds is essentially equivalent to the Bando-Siu conjecture \cite{BandoSiu:1994}, established by Daskalopoulos--Wentworth \cite{DaskWentworth:2004}, Sibling \cite{Sibley:2015} and Jacob \cite{Jacob:2013, Jacob:2013b}. In the context of $K$-stability and extremal Kähler metrics moment-weight inequalities are due to Tian \cite{Tian:1997}, Donaldson \cite{Donaldson:2002,Donaldson:2002b,Donaldson:2005} and Chen \cite{Chen:2008, Chen:2009}. In this context Chen--Sun \cite{ChenSun:2010} found an analytic proof of the Kempf existence theorem on finite dimensional spaces and we extend their argument to our infinite dimensional setting to prove existence of the dominant weight. The survey \cite{RobSaGeo} by Georgoulas--Robbin--Salamon provides an overview on the different proofs of the moment weight inequality for Hamiltoninan actions on closed Kähler manifolds and its importance for geometric invariant theory.

\subsubsection*{Acknowledgment}

I would like to thank my supervisor D. A. Salamon for many helpful discussions throughout the process of writing this paper.

\section{Preliminaries}

\subsection{The moment map picture}
We recall the natural Kähler structures on $\mathcal{A}(P)$ and $\mathcal{S}(P,X)$. Since $\mathcal{A}(P)$ is an affine space over the linear space $\Omega^1(\Sigma, \text{ad}(P))$, it suffices to specify the Kähler structure on the later one. For $a,b \in \Omega^1(\Sigma, \text{ad}(P))$ this is defined as
			$$\omega_{\mathcal{A}}(a,b) := \int_{\Sigma} \langle a \wedge b \rangle,\qquad J_{\mathcal{A}} a = *a, \qquad \langle a, b \rangle_{\mathcal{A}} := \int_{\Sigma} \langle a \wedge *b \rangle.$$
For $u \in \mathcal{S}(P,X)$ let $\tilde{u}: P \rightarrow X$ be the equivariant map determined by $u(z) = [p, \tilde{u}(p)]$ for $z \in \Sigma$ and $p \in P_z$. The tangent space $T_u \mathcal{S}(P,X)$ is represented by $G$-equivariant sections of the vector bundle $\tilde{u}^*TX \rightarrow P$ or equivalently by sections of the quotient bundle $\tilde{u}^*TX/G \rightarrow P/G = \Sigma$. The quotient bundle is again a vector bundle over $\Sigma$ and we denote it in the following by $u^*TX/G$ for simplicity. For $\hat{u}_1,\hat{u}_2 \in T_u \mathcal{S}(P,X) = \Omega^0(\Sigma, u^*TX/G)$ one defines
			$$\omega_{\mathcal{S}}(\hat{u}_1,\hat{u}_2) := \int_{\Sigma} \omega(\hat{u}_1, \hat{u}_2)\, dvol_{\Sigma},\quad J_{\mathcal{S}} \hat{u}_1 = J \hat{u}_1, \quad \langle \hat{u}_1, \hat{u}_1 \rangle_{\mathcal{S}} := \int_{\Sigma} \langle \hat{u}_1, \hat{u}_2 \rangle \, dvol_{\Sigma}.$$
On $\mathcal{A}(P)\times \mathcal{S}(P,X)$ denote the product Kähler structure by $(\omega_{\mathcal{A}\times\mathcal{S}}, J_{\mathcal{A}\times \mathcal{S}}, \langle \cdot , \cdot \rangle_{\mathcal{A}\times \mathcal{S}})$.

\begin{Lemma} \label{LemmaMoment}
The diagonal $\mathcal{G}(P)$-action on $\mathcal{A}(P)\times \mathcal{S}(P,X)$ is Hamiltonian with moment map
				\begin{align} \label{Momenteq0}\Phi: \mathcal{A}(P)\times\mathcal{S}(P,X) \rightarrow \Omega^0(\Sigma, \text{ad}(P)), \qquad \Phi(A,u) := *F_A + \mu(u). \end{align}
\end{Lemma}

\begin{proof} For $(A,u) \in \mathcal{A}(P)\times\mathcal{S}(P,X)$ and $\xi \in \Omega^0(\Sigma,\text{ad}(P))$ the infinitesimal action is given by
				\begin{align} \label{Momenteq1} 	\mathcal{L}_{(A,u)} \xi := \left.\frac{d}{dt}\right|_{t=0} \exp(t\xi) (A,u) = \left(-d_A \xi, L_u \xi\right) 	\end{align}		
where $L_x : \mathfrak{g}\rightarrow T_x X$ denotes the infinitesimal action of $\mathfrak{g}$ on $X$. The verification of (\ref{Momenteq1}) is straightforward and left to the reader. The differential of the function
					$$\mathcal{A}(P)\times\mathcal{S}(P,X) \rightarrow \mathbb{R}, \qquad (A,u) \mapsto \int_{\Sigma} \langle *F_A + \mu(u), \xi \rangle \, dvol_{\Sigma}$$
is the $1$-form 
				$$T_A \mathcal{A}(P) \times T_u \mathcal{S}(P,X) \rightarrow \mathbb{R}, \qquad
										(a,\hat{u}) \mapsto \int_{\Sigma} \langle - d_A\xi \wedge a \rangle + \int_\Sigma \omega( L_u\xi, \hat{u}).$$
This is precisely $\omega_{\mathcal{A}\times\mathcal{S}}(\mathcal{L}_{(A,u)}\xi, \cdot )$ and (\ref{Momenteq0}) satisfies the moment map equation.
\end{proof}

\subsection{Connections on \texorpdfstring{$P(X)$}{PX} and the space \texorpdfstring{$\mathcal{S}(P,X)$}{SPX}}

For $x \in X$ the infinitesimal action of $\mathfrak{g}$ defines a map $L_x : \mathfrak{g} \rightarrow T_x X$. A smooth connection $A \in \mathcal{A}(P)$ induces on the Kähler fibration $P(X)$ the covariant derivative
			$$d_A: \Omega^0(\Sigma, P(X)) \rightarrow \Omega^1(\Sigma, u^*TX/G), \qquad d_A u := d u + L_u A$$
with values in the the vertical tangent bundle along $u$ which is isomorphic to $u^*TX/G$. Moreover, $A$ and the Levi-Civita connection induce a covariant derivative
			$$\nabla_A :\Omega^0(\Sigma, u^*TX/G) \rightarrow \Omega^1(\Sigma, u^*TX/G), \qquad \nabla_A \xi := \nabla \xi  + \nabla_{\xi} (L_u A).$$
All these covariant derivatives extend to first order elliptic operators between suitable Sobolev spaces.

\subsection{The Yang--Mills--Higgs functional}

The moment map squared functional (\ref{DefF}) and the Yang--Mills--Higgs functional (\ref{DefYMH}) are related by the following energy identity.

\begin{Proposition}\label{PropEngergy}
Let $(A,u) \in \mathcal{A}(P)\times \mathcal{S}(P,X)$, then
		\begin{align} \label{energyid}
				\mathcal{YMH}(A,u) = \mathcal{F}(A,u) + \int_{\Sigma} ||\bar{\partial}_A(u)||^2\, dvol_{\Sigma} + \langle \omega - \mu, u \rangle
		\end{align}
where	$\langle \omega - \mu, u \rangle = \int_{\Sigma} u^*\omega - d \langle \mu(u), A \rangle$.
\end{Proposition}

\begin{proof} This is Proposition 3.1 in \cite{Salamon:2000}. \end{proof}

\begin{Remark} The term $\langle \omega - \mu, u \rangle$ describes the pairing between the equivariant homology class $[u] \in H_2^G(X, \mathbb{Z})$ determined by $u$ and the equivariant cohomology class $[\omega - \mu] \in H^2_G(X, \mathbb{R})$, see \cite{Salamon:2000} for more details. In particular, this term is constant on the homotopy class of $(A, u)$ and solutions of the symplectic vortex equation (\ref{vortexeq1}) minimize the Yang--Mills--Higgs functional in their homotopy class.
\end{Remark}

\begin{Lemma} \label{LemmaGrad}
	
	\begin{enumerate}
			\item The $L^2$-gradient of $\mathcal{F}$ is given by
										$$\nabla \mathcal{F}(A,u) = \begin{pmatrix} -*d_A( *F_A + \mu(u)) \\ J L_u(*F_A + \mu(u)) \end{pmatrix}.$$
			\item The $L^2$-gradient of $\mathcal{YMH}$ is given by 
										$$\nabla \mathcal{YMH}(A,u) = \begin{pmatrix} d_A^*F_A + L_u^* d_A u \\ \nabla_A^* d_A u + d\mu(u)^* \mu(u) \end{pmatrix}.$$
			\item If $(A,u) \in \mathcal{H}(P,X)$, then both gradients are tangential to $\mathcal{H}(P,X)$ and agree. That is
					\begin{align*} -*d_A( *F_A + \mu(u))   &= d_A^*F_A + L_u^* d_A u   \\
													J(u) L_u(*F_A + \mu(u))	&=	 \nabla_A^* d_A u + d\mu(u)^* \mu(u)
					\end{align*}
			holds for all $(A,u) \in \mathcal{H}(P,X)$.
	\end{enumerate}		
\end{Lemma}

\begin{proof}
We leave the first two parts to the reader (or refer to \cite{Venugopalan:2016} and \cite{Lin:2012} for full details). For the last claim, note that
			$$\nabla \mathcal{F}(A,u) = \begin{pmatrix} -*d_A( *F_A + \mu(u)) \\ J(u) L_u(*F_A + \mu(u)) \end{pmatrix} = \mathcal{L}_{(A,u)}^c \textbf{i}(*F_A + \mu(u))$$
is tangential to the complexified orbit $\mathcal{G}^c (A,u) \subset \mathcal{H}(P,X)$ and hence tangential to $\mathcal{H}(P,X)$. Since $\mathcal{H}(P,X)$ minimizes the functional 
				$$(A,u) \mapsto \int_{\Sigma} ||\bar{\partial}_A u||^2\, dvol_{\Sigma}$$
its gradient vanishes for $(A,u) \in \mathcal{H}(P,X)$ and the claim follows from the energy identity (\ref{energyid}).
\end{proof}

\subsection{Equivariant convexity at infinity}

The next Lemma shows that under assumption \textbf{(C)} solutions of (\ref{floweq0}) remain in a compact region of $X$.

\begin{Lemma}\label{LemmaEC}
Suppose $X$ is equivariantly convex at infinity and let $f: X \rightarrow [0,\infty)$ and $c_0 > 0$ be as in (\ref{EqC}). Let $T \in (0,\infty]$ and suppose $(A,u) : [0,T) \rightarrow \mathcal{H}(P,X)$ is a smooth map satisfying
						$$\partial_t u = -JL_u(*F_A + \mu(u)).$$
Then, for $c > c_0$ and $S_c := f^{-1}[0,c]$, it holds
			$$u_0(P) \subset S_c \qquad \Longrightarrow \qquad u_t(P) \subset S_c$$
for every $t \in [0,T]$.					
\end{Lemma}

\begin{proof}
The proof is similar to the calculation in \cite{Salamon:2002}, Lemma 2.7.

In local trivializing coordinates $z = x + \textbf{i} y$ define
			$$v_x := \partial^A_x u:= \partial_x u + L_u A(\partial_x), \qquad v_y := \partial^A_y u:= \partial_y u + L_u A(\partial_y)$$
Denote by $\tilde{\Delta} := \partial_x^2 + \partial_y^2$ the standard Laplacian. Then
		\begin{align*}
				\tilde{\Delta} f(u) &= \partial_x \langle \nabla f(u), v_x \rangle + \partial_y \langle \nabla f(u), v_y \rangle \\
										&= \langle \nabla^A_x \nabla f(u), v_x \rangle + \langle  \nabla^A_y \nabla f(u), v_y \rangle + \langle \nabla f(u), \nabla^A_x v_x + \nabla^A_y v_y \rangle
		\end{align*}
and since $f$ is $G$-invariant, we obtain
		\begin{align} \label{ECeq1}
				\tilde{\Delta} f(u) &= \langle \nabla_{v_x} \nabla f(u), v_x \rangle + \langle \nabla_{v_y} \nabla f(u), v_y \rangle + \langle \nabla f(u), \nabla^A_x v_x + \nabla^A_y v_y \rangle
		\end{align}
Using the characteristic equation for the curvature
			\begin{align*}  \nabla^A_x v_y - \nabla^A_y v_x = L_u F_A(\partial_x,\partial_y) \end{align*}
and the assumption $(A,u) \in \mathcal{H}(P,X)$, which is equivalent to $v_x + J v_y = 0$, we obtain
			$$\nabla^A_x v_x + \nabla^A_y v_y = -J \left(\nabla^A_x v_x - \nabla^A_y v_x \right) = - J L_u F_A(\partial_x, \partial_y).$$
Inserting this in	(\ref{ECeq1}) yields
		\begin{align} \label{ECeq2}
				\tilde{\Delta} f(u) &= \langle \nabla_{v_x} \nabla f(u), v_x \rangle + \langle \nabla_{v_y} \nabla f(u), v_y \rangle - \langle \nabla f(u), JL_u F_A(\partial_x, \partial_y) \rangle
		\end{align}	
If $f(u) \geq c_0$, then the convexity assumption implies that the first two terms in (\ref{ECeq2}) are positive and thus
			$$f(u) \geq c_0 \qquad \Longrightarrow \qquad \Delta f(u) \leq \langle \nabla f(u), J L_u *F_A \rangle $$
where $\Delta = d^*d$ denotes the positive Laplacian (which corresponds to $- \lambda \tilde{\Delta}$ in local coordinates for some function $\lambda > 0$). This yields
		\begin{align} \label{ECeq3} f(u) \geq c_0 \qquad \Longrightarrow \qquad \left(\partial_t + \Delta \right) f(u) \leq \langle \nabla f(u), - J L_u d\mu(u) \rangle \leq 0 \end{align}
where we used the second equation in the convexity assumption. 

We deduce the claim from (\ref{ECeq3}) by contradiction. Suppose there exists $M > c$ such that
					$$t_1 := \inf \{t \in [0,T) \,|\, \text{$f(u_t(z)) \geq M$ for some $z \in \Sigma$}\}$$
satisfies $0 < t_1 < T$ (i.e. $\inf \emptyset = \infty$ is excluded). Let $D \subset \Sigma$ be a small disc and let $t_0 \in (0,t_1)$ be such that 
					$$f(u_t(z)) > c_0\qquad \forall (t,z)  \in [t_0, t_1]\times D$$
and $f(u_{t_1}(z_0)) = M$ for some interior point $z_0 \in D$. It follows from (\ref{ECeq3}) that in local coordinates $f(u_t(z))$ is a subsolution to a parabolic equation on $[t_0,t_1]\times D$ and by construction it attains its maximum on $\{t_1\}\times D$. By the strong maximum principle for parabolic equations (see \cite{Friedman:ParabolicPDE} Chapter 2, Theorem 1), it follows that $f(u_t(x)) \equiv M$ is constant on $[t_0,t_1]\times D$. This contradicts the definition of $t_1$ and completes the proof of the Lemma.

\end{proof}

\subsection{Sobolev spaces}

We discuss mixed Sobolev spaces of time dependent sections of vector bundles. Following R\r{a}de \cite{Rade1992} and Venugopalan \cite{Venugopalan:2016} we shall use a norm on $H^{r}([0,t_0], H^s(\Sigma,V))$ which depends on the length $t_0$ of the time interval. For convenience, we use the abbreviation $H^{s} := W^{s,2}$ for $L^2$-Sobolev spaces.

\subsubsection{Fractional Sobolev spaces on bounded domains}
The refer to \cite{Adams:Sobolev} for the general theory of Sobolev spaces. The definition of fractional Sobolev spaces (also called Bessel potential spaces) uses deep results from harmonic analysis (see \cite{Stein:SingInt} Chapter V.3 or \cite{Hamilton:HeatMonograph} Chapter 2.1-3). For $s \in \mathbb{R}$ and $p \in (1,\infty)$ one defines
			\begin{align} 
					\label{SobolevDefn1} W^{s,p}(\mathbb{R}^n)  := (1 - \Delta)^{-s/2} \left(L^p(\mathbb{R}^n) \right), \qquad ||f||_{W^{s,p}} := ||(1 - \Delta)^{s/2} f ||_{L^p}.
			\end{align}
For a smoothly bounded domain $\Omega \subset \mathbb{R}^n$ and $f \in C^{\infty}(\Omega)$ one defines
				$$||f||_{W^{s,p}(\Omega)} = \inf_{f= F|_\Omega} ||F||_{W^{s,p}(\mathbb{R}^n)}$$
where the infimum ranges over all $F \in C^{\infty}_0(\mathbb{R}^n)$ which restrict to $f$. The space $W^{s,p}(\Omega)$ (resp. $W^{s,p}_0(\Omega)$) is the closure of $C^{\infty}(\Omega)$ (resp. $C^{\infty}_0(\Omega)$) under this norm. The extension theorem shows that $W^{s,p}(\Omega)$ is the set of restriction to $\Omega$ of functions in $W^{s,p}(\mathbb{R}^n)$. In the special case $p = 2$ one obtains the Hilbert spaces $H^s(\Omega) = W^{s,2}(\Omega)$ (see \cite{LionsMangenes1}).

\paragraph{Interpolation.} The spaces $W^{s,p}(\Omega)$ form a family of interpolation spaces in both parameters: the degree $s$ of differentiability and the degree $p$ of summability. For $1 < p_0, p_1 < \infty$,  $s_0, s_1 \in \mathbb{R}$ and $0 < \theta < 1$ it holds
		\begin{align} \label{SobolevEq1}
				W^{s_{\theta},p_{\theta}}(\Omega) \cong  [W^{s_0,p_0}(\Omega), W^{s_1,p_1}(\Omega)]_{\theta} 
		\end{align}		
with $s_{\theta} = (1- \theta) s_0 + \theta s_1$, $p_{\theta} = (1- \theta) p_0 + \theta p_1$ and $[\cdot, \cdot]_{\theta}$ refers to the holomorphic interpolation method. The same remains valid for the spaces $W^{s,p}_0(\Omega)$. (See \cite{Hamilton:HeatMonograph} Chapter 2.4-5, \cite{Triebel:Interpol} Chapters 1.9, 2.4 and 4.3).

\paragraph{Duality.} For $1 < p,q < \infty$, $\frac{1}{p} + \frac{1}{q} = 1$ and $s \in \mathbb{R}^+_0$ there exists a natural identification
		\begin{align} \label{SobolevEq2}
				W^{-s,p}(\Omega) \cong W^{s, q}_0(\Omega)^*.
		\end{align}
which is obtained by extending the $L^2$-product.

\paragraph{Products.} Let $s,t,u \in \mathbb{R}$ with $u \leq \min\{s,t\}$ and $s +t \geq 0$. Let $1 < p,q,r < \infty$ with $s \neq n/p$, $t \neq n/q$, $u \neq -n/r$, and
			\begin{align} \label{pcond}
					\max\left\{ \left(\frac{1}{p} - \frac{s}{n}\right), \, \left(\frac{1}{q} - \frac{t}{n}\right),\, \left(\frac{1}{p} - \frac{s}{n}\right) + \left(\frac{1}{q} - \frac{t}{n}\right) \right\} \leq \left(\frac{1}{r} -\frac{u}{n}\right)
			\end{align}		
Then, if $f \in W^{s,p}(\Omega)$ and $g \in W^{t,q}(\Omega)$, the product $fg$ is contained in $W^{u,r}(\Omega)$ and satisfies an estimate
			\begin{align} \label{peq}
					||fg||_{W^{u,r}(\Omega)} \leq C ||f||_{W^{s,p}(\Omega)} ||g||_{W^{t,q}(\Omega)}.
			\end{align}		
This follows for $s,t,u \in \mathbb{Z}^+_0$ from the Sobolev embedding theorem and Hölder's inequality. The general case is obtained from this by interpolation (\ref{SobolevEq1}) and duality (\ref{SobolevEq2}). (See \cite{Palais:1968} Theorem 9.6 for the details)

\subsubsection{Sobolev spaces of sections \texorpdfstring{$H^s(\Sigma, V)$}{HsV}}

Let $V \rightarrow \Sigma$ be a Riemannian vector bundle over $\Sigma$. One can describe $H^s(\Sigma, V)$ in local coordinates as follows: Let $\{U_{\alpha}\}$ be an open trivializing cover of $\Sigma$ by charts and choose unitary trivializations $V|_{U_{\alpha}} \cong U_{\alpha} \times \mathbb{R}^n$. A partition of unity subordinate to the cover $\{U_{\alpha}\}$ divides a section $\sigma \in \Omega^0(\Sigma, V)$ into a collection of functions $\sigma_{\alpha}^j \in C^{\infty}_0(U_{\alpha})$. Using the charts we identify $U_{\alpha}$ with open bounded subsets $\Omega_{\alpha} \subset \mathbb{R}^2$ and define
					\begin{align} \label{SobolevDef2} ||\sigma||_{H^{s}} := \sum_{\alpha, j} ||\sigma_{\alpha}^j||_{H^{s}(\Omega_{\alpha})}. \end{align}
The space $H^{s}(\Sigma, V)$ is the completion of $\Omega^0(\Sigma, V)$ in this norm.

\begin{Remark}
Let $\nabla$ be a smooth Riemannian connection on $V$. For $s = k \in \mathbb{Z}^+_0$ the norm
				\begin{align} \label{SobolevDef3} ||\sigma|| := \sum_{j=0}^k ||\nabla^j \sigma||_{L^2} \end{align}
is equivalent to the $H^{k}$-norm defined in (\ref{SobolevDef2}). This leads to an alternative construction of $H^{s}(\Sigma, V)$ starting with (\ref{SobolevDef3}) for positive integers and then using interpolation and duality.
\end{Remark}

The product formula (\ref{peq}) takes under the assumptions $p=q=r=2$ and $n=2$ the following simpler form. Let $s,t,u \in \mathbb{R}$ with $s,t \neq +1$, $u \neq -1$, $s + t \geq 0$, and
							$$u \leq \min\{s,t, s+t - 1\}.$$
Then, if $f \in H^s(\Sigma)$ and $g \in H^t(\Sigma)$, the product $fg$ is contained in $H^u(\Sigma)$ and satisfies an estimate
							$$||fg||_{H^u} \leq C ||f||_{H^s} ||g||_{H^t}.$$

\subsubsection{Time dependent Sobolev spaces \texorpdfstring{$H^{r}([0,t_0],H^s(\Sigma,V))$}{HrHsV}}
Let $K$ be a separable Hilbert space, let $t_0 > 0$ and let $f: [0,t_0] \rightarrow K$ be a smooth function. The following is a slight variant of (\ref{SobolevDefn1}). For $r \in \mathbb{R}$ we define
					\begin{align} \label{SobolevDefn3}
							||f||_{H^r([0,t_0], K)} := \inf_{F|_{[0,t_0]} = f} \left( \int_{-\infty}^\infty (\tau^2 + t_0^{-2})^{r} ||\hat{F}(\tau)||_K^2 \,d\tau \right)^{\frac{1}{2}}
					\end{align}		
where the infimum is taken over all $F \in C_0^{\infty}(\mathbb{R}, K)$ which restrict to $f$ on $[0,t_0]$ and $\hat{F}$ denotes the Fourier transform.

\begin{Remark}
If $r = k \in \mathbb{Z}^+_0$, then (\ref{SobolevDefn3}) is equivalent to the norm
					\begin{align} \label{SobolevDefn4}
							||f|| := \sum_{j=0}^k\left|\left|t_0^{-(k-j)} \frac{d^j}{dt^j} f \right|\right|_{L^2([0,t_0],K)}^2.
					\end{align}		
As before, one could construct the spaces $H^r([0,t_0], K)$ using (\ref{SobolevDefn4}) for positive integers and then use interpolation and duality.
\end{Remark}

The dependence of the norms on $t_0$ has the advantage that for $r_1 \geq r_2$ the inclusion
				$$H^{r_1}([0,t_0], K) \hookrightarrow H^{r_2}([0,t_0], K)$$
has norm $\leq C t_0^{r_1 - r_2}$. In particular, this norm can be controlled by $t_0$.

For $K = H^s(\Sigma, V)$ one obtains the spaces $H^{r}([0,t_0],H^s(\Sigma,V))$. These form again a family of interpolation spaces (\cite{Venugopalan:2016} Lemma 6.36): For $s_0, s_1, r_0, r_1 \in \mathbb{R}$ and $\theta \in (0,1)$ we have
			$$\left[H^{r_0}([0,t_0],H^{s_0}(\Sigma, V)), H^{r_1}([0,t_0],H^{s_1}(\Sigma, V)) \right]_{\theta} \cong H^{r_{\theta}}([0,t_0],H^{s_{\theta}}(\Sigma,V))$$
with $r_{\theta} = (1-\theta)r_0 + \theta r_1$, $s_{\theta} = (1-\theta)s_0 + \theta s_1$ and $[\cdot,\cdot]_{\theta}$ denotes the holomorphic interpolation method.

\subsection{The heat equation}
Let $V \rightarrow \Sigma$ be a Riemannian vector bundle and let $\nabla$ be a Riemannian connection on $V$.
		
\begin{Lemma} \label{LemmaHeat1}
For every $\sigma_0 \in \Omega^0(\Sigma, V)$ and $t_0 > 0$, there exists a unique smooth solution $\sigma: [0,t_0] \rightarrow \Omega^0(\Sigma, V)$ solving the initial value problem
		\begin{align} \label{heateq1}
						\partial_t \sigma + \nabla^* \nabla \sigma = 0, \qquad \sigma(0,\cdot) = \sigma_0.
		\end{align}
Moreover, there exists a constant $C > 0$	such that the following estimate holds
					$$||\sigma||_{L^2([0,t_0], H^1(\Sigma, V))} \leq C t_0^{\frac{1}{2}} ||\sigma_0||_{L^2}.$$
\end{Lemma}

\begin{proof}
This is a special case of Lemma 6.33 in \cite{Venugopalan:2016}.
\end{proof}

From this we deduce the following estimates.

\begin{Lemma} \label{LemmaHeat2}
Let $f :[0,t_0] \rightarrow \Omega^0(\Sigma, V)$ be smooth. There exists a unique smooth solution $\psi$ of the equation
				\begin{align} \label{heateq2}
							\partial_t \psi + \nabla^* \nabla \psi = f, \qquad \psi(0,\cdot) = 0.
				\end{align}
Moreover, the solution satisfies the estimates
			\begin{align} \label{heateq3} ||\psi||_{L^2([0,t_0], H^1(\Sigma, V))} \leq C t_0^{\frac{1}{2}}||\psi||_{L^1([0,t_0], L^2(\Sigma, V))} \end{align}
and
			\begin{align} \label{heateq4} ||\psi||_{L^2([0,t_0], H^1(\Sigma, V))} \leq C t_0^{\frac{1}{4}} ||\psi||_{L^2([0,t_0], H^{-\frac{1}{2}}(\Sigma, V))} .\end{align}
\end{Lemma}

\begin{proof}
Let $P_t$ denote the solution operator of (\ref{heateq1}), i.e. $P_0 = \mathds{1}$ and $P_t \sigma_0(\cdot) = \sigma(t,\cdot)$ satisfies (\ref{heateq1}). The solution of (\ref{heateq2}) is then given by
						$$\psi(t,\cdot) = \int_0^t P_{t-s} f(s,\cdot)\, ds.$$
The Minkowski inequality and Lemma \ref{LemmaHeat1} yield
		\begin{align*}
			||\psi||_{L^2(H^1)} &\leq \left( \int_0^{t_0} \left(\int_0^{t} ||P_{t-s} f(s,\cdot)||_{H^1}\, ds \right)^2 dt \right)^{\frac{1}{2}} \\
													&\leq \int_0^{t_0} \left(\int_s^{t_0} ||P_{t-s} f(s,\cdot)||_{H^1}^2\, dt \right)^{\frac{1}{2}} ds \\
													&\leq C t_0^{\frac{1}{2}}\int_0^{t_0} ||f(s,\cdot)||_{L^2}\, ds
		\end{align*}
and this proves (\ref{heateq3}). 
Abbreviate 
			$$H^{r,s} := H^r([0,t_0], H^s(\Sigma, V)).$$ 
Parabolic regularity (see \cite{Venugopalan:2016} Lemma 6.35) yields
			$$||\psi||_{H^{\frac{3}{4}, -\frac{1}{2}}} \leq C||f||_{H^{-\frac{1}{4}, -\frac{1}{2}}}, \qquad ||\psi||_{H^{-\frac{1}{4},\frac{3}{2}}} \leq C||f||_{H^{-\frac{1}{4},-\frac{1}{2}}}$$ 
and, since $H^{0,1}$ is an interpolation space between $H^{\frac{3}{4},-\frac{1}{2}}$ and $H^{-\frac{1}{4},\frac{3}{2}}$, it follows
			$$||\psi||_{H^{0,1}} \leq C ||f||_{H^{-\frac{1}{4},-\frac{1}{2}}} \leq C t_0^{\frac{1}{4}} ||f||_{H^{0,-\frac{1}{2}}}.$$
This establishes (\ref{heateq4}) and completes the proof.
\end{proof}

\subsection{Sobolev completions and regularity assumptions} \label{Section.SobolevReg}
For the main part of the article, we need to consider suitable Sobolev completions of the various spaces defined in the introduction. The space 
				$$\mathcal{S}^{2,2}(P,X) := W^{2,2}(\Sigma, P(X))$$ 
contains all continuous sections $u: \Sigma \rightarrow P(X)$ which in any trivialization of $P(X)$ and local coordinates on $\Sigma$ and $X$ are of Sobolev class $W^{2,2}_{loc}$. It carries a natural topology, since $W^{2,2}_{loc}(\mathbb{R}^2) \hookrightarrow C^0(\mathbb{R}^2)$ is in the good range of the Sobolev embedding: For $u \in \mathcal{S}(P,X)$ let $\epsilon > 0$ be smaller then the injectivity radius of $X$ along the image of $u$. Then
				$$\{ \hat{u} \in W^{2,2}(\Sigma, u^* TX/G)\, |\, ||\hat{u}||_{W^{2,2}} < \epsilon \} \rightarrow \mathcal{S}^{2,2}(P,X), \qquad \hat{u} \mapsto \exp_{u} \hat{u}$$
defines a homeomorphism onto its image. 

With respect to a smooth reference connection $A_0 \in \mathcal{A}(P)$, we define
			$$\mathcal{A}^{1,2}(P) := \{ A_0 + a \,|\, a \in W^{1,2}(\Sigma, T^*\Sigma \otimes \text{ad}(P))\}$$
and denote
			$$\mathcal{H}^{1,2}(P,X) := \{ (A,u) \in \mathcal{A}^{1,2}(P)\times \mathcal{S}^{2,2}(P,X)\,|\, \bar{\partial}_A u = 0 \}.$$
The $W^{2,2}$ completion of the gauge groups
			$$\mathcal{G}^{2,2}(P) := W^{2,2}(\Sigma, \text{Ad}(P)), \qquad (\mathcal{G}^c)^{2,2}(P) := W^{2,2}(\Sigma, \text{Ad}(P\times_G G^c))$$
are similar defined as $\mathcal{S}^{2,2}(P,X)$ by requiring their sections to be of Sobolev class $W^{2,2}$ in any local trivialization. These groups act continuously on $\mathcal{S}^{2,2}(P,X)$, $\mathcal{A}^{1,2}(P)$ and $\mathcal{H}^{1,2}(P,X)$ as one readily checks.

\begin{Lemma}\label{Lemma.SobolevReg}
Let $(A,u) \in \mathcal{H}^{1,2}(P,X)$.
	\begin{enumerate}	
			\item There exists $g \in (\mathcal{G}^c)^{2,2}(P)$ such that $g(A,u)$ is smooth.
			\item If $(A,u)$ is a critical point of $\mathcal{YMH}$ satisfying
							\begin{align} \label{CritEq} d_A(*F_A + \mu) = 0, \qquad L_u(*F_A + \mu(u)) = 0 \end{align}
				then there exits $k \in \mathcal{G}^{2,2}(P)$ such that $k(A,u)$ is smooth.
	\end{enumerate}
\end{Lemma}

\begin{proof}
This Lemma is proven as in the Yang-Mills case. First, there exists $g \in (\mathcal{G}^c)^{2,2}$ such that $gA$ is smooth (see \cite{AtBott:YangMillsEq}, Lemma 14.8). Then $\bar{\partial}_{gA} (gu) = 0$ and elliptic regularity yields that $gu$ is smooth. This proves the first part of the Lemma.

For the second part we pass to a Coulomb gauge and choose a smooth reference connection $A_0 \in \mathcal{A}(P)$ and $k \in \mathcal{G}^{2,2}(P)$ such that $d_{A_0}^*(kA - A_0) = 0$. By (\ref{CritEq}), $a := kA - A_0$ satisfies
			\begin{align} \label{bootstrap} \Delta_{A_0} a = d_{A_0}* F_{A_0} + \frac{1}{2} [a\wedge a] + d_{A_0}(\mu(ku)) + [a \wedge (*F_{kA} + \mu(ku))]. \end{align}
Suppose first that $a \in H^1$ and $u \in H^{2}$. Using the multiplication theorem $H^{1}\otimes L^2 \rightarrow H^{-\frac{1}{2}}$, one sees that the right hand side of (\ref{bootstrap}) is in $H^{-\frac{1}{2}}$ and hence $a \in H^{\frac{3}{2}}$. With this improved regularity, the right hand side of (\ref{bootstrap}) is now contained in $H^0$ and hence $a \in H^2$. The holomorphicity condition
				$$0 = \bar{\partial}_{kA} (ku) = \bar{\partial}_{A_0}(ku) + \left( L_{ku} a \right)^{0,1}$$
then yields $u \in H^3$. Repeating this argument, one shows $k(A,u) \in H^{\ell}\times H^{\ell +1}$ for every $\ell \geq 2$ and this completes the bootstrapping argument.
\end{proof}

\section{Convergence of the Yang--Mills--Higgs flow}

In the first section weak solutions of the gradient flow (\ref{floweq0}) are defined and the existence and regularity of solutions are discussed. The second section contains a proof of the {\L}ojasiewicz gradient inequality for the Yang--Mills--Higgs functional. Combining this inequality with an interior regularity result in the third section, we can then prove that solutions convergence under the additional assumptions \textbf{(A)}, \textbf{(B)}. This approach is very similar to the one developed by R\r{a}de \cite{Rade1992} in the Yang--Mills case.

\subsection{The gradient flow equations}

\begin{Definition}[\textbf{Negative gradient flow of $\mathcal{YMH}$}] A (weak) solution of
			\begin{align} \label{floweq1}
					\begin{split}	
										 \partial_t A &= -d_A^*F_A - L_u^* d_A u \\
										 \partial_t u &= - \nabla_A^* d_A u - d\mu(u)^* \mu(u)
					\end{split}					
			\end{align}
is a continuous map $(A,u) : [0,\infty) \rightarrow \mathcal{H}^{1,2}(P,X)$, such that there exists a sequence of smooth solutions of (\ref{floweq1}) converging to $(A,u)$ in $C^0([0,\infty), H^1\times H^2)$.
\end{Definition}

\begin{Definition}[\textbf{Negative gradient flow of $\mathcal{F}$}] A (weak) solution of
			\begin{align} \label{floweq2}
					\begin{split}	
										 \partial_t A &= *d_A( *F_A + \mu(u)) \\
										 \partial_t u &= - J L_u( *F_A + \mu(u))
					\end{split}					
			\end{align}										
is a continuous map $(A,u) : [0,\infty) \rightarrow \mathcal{H}^{1,2}(P,X)$, such that there exists a sequence of smooth solutions of (\ref{floweq2}) converging to $(A,u)$ in $C^0([0,\infty), H^1\times H^2)$.
\end{Definition}

By Lemma \ref{LemmaGrad} both of these flows agree:
			$$\text{$(A,u)$ is a weak solution of (\ref{floweq1})} \quad \Longleftrightarrow \quad \text{$(A,u)$ is a weak solution of (\ref{floweq2})}.$$ 
The following theorem is a slight extension of a result of Venugopalan \cite{Venugopalan:2016} (she works in the $H^1\times C^0$ topology and needs to assume that the flow remains in a compact region of $X$).

\begin{Theorem} \label{ThmVen}
Assume \textbf{(C)} and let $(A_0,u_0) \in \mathcal{H}^{1,2}(P,X)$.
	\begin{enumerate}	
		\item There exists a unique solution $(A,u) \in C^0([0,\infty), \mathcal{H}^{1,2}(P,X))$ of (\ref{floweq2}) with $A(0) = A_0$ and $u(0,\cdot) = u_0$.
		\item The map $\Phi: [0,\infty) \rightarrow L^2(\Sigma, T^*\Sigma\otimes \text{ad}(P))$
								\begin{align} \label{eqV1} \Phi(t) := *F_{A(t)} + \mu(u(t)) \end{align}
					is contained in the spaces $C^0([0,\infty), L^2)$ and $L^2_{loc}([0,\infty), H^1)$.
		\item The solution $g: [0,\infty) \rightarrow (\mathcal{G}^c)^{2,2}(P)$ of the ODE
									\begin{align}	\label{KNflow} g^{-1}(t)\dot{g}(t) = \textbf{i} (*F_{A(t)} + \mu(u(t))), \qquad g(0) = \mathds{1}. \end{align}
					is continuous with values in $H^2$ and satisfies $(A(t),u(t)) = g(t)^{-1}(A_0, u_0)$.
		\item The solution $(A(t),u(t))$ of (\ref{floweq2}), the map $\Phi(t)$ in (\ref{eqV1}) and the solution $g(t)$ of (\ref{KNflow}) depend continuously on the initial condition $(A_0,u_0) \in \mathcal{H}^{1,2}(P,X)$ in the respective topologies stated above.
	\end{enumerate}
\end{Theorem}

\begin{proof}
Let $f: X \rightarrow [0,\infty)$ and $c_0 \in \mathbb{R}$ be as in (\ref{EqC}). By Lemma \ref{LemmaEC} the compact sets $S_c := f^{-1}[0,c]$ with $c > c_0$ have the following property: If $(A(t),u(t))$ is a gradient flow line starting at $(A_0, u_0)$, then
			\begin{align} \label{CompEq1} u_0(P) \subset S_c \qquad \Longrightarrow \qquad u_t(P) \subset S_c \quad \forall t \geq 0. \end{align}
Venugopalan proves long time existence by establishing short time existence together with an uniform lower bound on the existence intervall. When we restrict to the set $S_c$ her analysis yields uniform lower bounds for the existence interval for any solution with $u_0(P) \subset S_c$. Now Theorem 1.1 in \cite{Venugopalan:2016} shows that for any $A_0 \in H^1$ and $u_0 \in C^0$ there exists a unique (weak) solution $(A,u) \in C^0([0,\infty), H^1\times C^0)$. Moreover, the proof shows the solution $(A,u)$ depends continuously on the initial condition $(A_0,u_0)$, the moment map term $\Phi(t) := *F_{A(t)} + \mu(u(t))$ is contained in the space $C^0([0,\infty), L^2) \cap L^2_{loc}([0,\infty), H^1)$ and depends continously on the initial condition $(A_0,u_0)$ in these topologies. The additional regularity $\Phi \in L^2_{loc}([0,\infty), H^1)$ is somewhat hidden in her proof and follows from the consideration of the space $\tilde{U}_P(t_0)$ at the end of the proof of Proposition 3.3. There she shows $\Phi \in H^{\frac{1}{2} + \epsilon, -2\epsilon}\cap H^{-\frac{1}{2}, 2}$ and this embedds into $L^2(H^1)$ by interpolation.

By the Sobolev embedding $H^1\times H^2 \hookrightarrow H^1 \times C^0$, we obtain for any initial condition $(A_0,u_0) \in H^1 \times H^2$ a solution $(A,u) \in C^0([0,\infty), H^1\times C^0)$. We claim that there exists a continuous path of complex gauge transformations $g : [0,\infty) \rightarrow (\mathcal{G}^c)^{2,2}(P)$, depending continuously on the initial condition $(A_0,u_0)$, such that $(A(t), u(t)) = g(t)^{-1}(A_0,u_0)$. By continuity of the gauge action, this readily implies that $(A,u) \in C^0([0,\infty), H^1\times H^2)$ and it depends continuously on the initial condition.

For the claim, $\Phi \in L^2_{loc}([0,\infty), H^1)$ and (\ref{KNflow}) yield $g \in H^1_{loc}([0,\infty), H^1)$. Moreover (\ref{floweq2}) shows $\partial_t A \in L^2_{loc}([0,\infty), L^2)$. Hence $B(t) := A(t) - g(t)^{-1}A_0 \in H^1_{loc}([0,\infty), L^1)$ and
								$$\dot{B}(t) = [*B(t), \Phi(t)], \qquad B(0) = 0.$$
If $B$ is smooth, this implies $B = 0$. In general, one can approximate weak solutions by smooth solutions and deduce then $B = 0$. This shows $A(t) = g(t)^{-1}A_0$ for all $t \geq 0$. Since $g \in H^1_{loc}([0,\infty), H^1)$ and $A \in C^0([0,\infty), H^1)$ depend continuously on the initial condition in $H^1\times H^2$, it follows from the equation
			$$A(t)^{0,1} = (g(t)^{-1}A_0)^{0,1} = A_0 + g(t)^{-1} \bar{\partial}_{A_0} g(t)$$
and standard elliptic bootstrapping arguments that $g \in C^0([0,\infty), H^2)$ depends continuously on the initial condition. One readily checks $u(t) = g(t)^{-1} u_0$ and this completes the proof.
\end{proof}

\subsection{{\L}ojasiewicz gradient inequality}

The {\L}ojasiewicz gradient inequality is the key ingredient in proving uniform convergence of the Yang--Mills--Higgs flow. This approach is due to Simon \cite{Simon:1983} and we follow quite closely the arguments of R\r{a}de \cite{Rade1992} in the Yang--Mills case.

\begin{Theorem}[\textbf{{\L}ojasiewicz gradient inequality}] \label{ThmLoj1}
Assume \textbf{(A)} and let $(A_{\infty},u_{\infty}) \in \mathcal{H}^{1,2}(P,X)$ be a critical point of $\mathcal{YMH}$. Then there exist $\epsilon, C > 0$ and $\gamma \in [\frac{1}{2}, 1)$ such that for all $a \in H^1(\Sigma, T^*\Sigma \otimes \text{ad}(P))$ and $\hat{u} \in H^2(\Sigma, u^*TX/G)$ with  $||a||_{H^1} + ||\hat{u}||_{H^2} < \epsilon$ it holds
			\begin{equation} \label{Lojeq1}
						\begin{aligned}	&||\nabla \mathcal{YMH}(A_{\infty} + a, \exp_{u_{\infty}} \hat{u})||_{H^{-1}\times L^2} \\
														&\qquad \geq c | \mathcal{YMH}(A_{\infty} + a, \exp_{u_{\infty}} \hat{u}) - \mathcal{YMH}(A_{\infty},u_{\infty}) |^{\gamma}.
						\end{aligned}
			\end{equation}
\end{Theorem}

\begin{proof}
See page \pageref{proofThmLoj1}.
\end{proof}

By Lemma \ref{Lemma.SobolevReg} every critical point $(A_0,u_0) \in \mathcal{H}^{1,2}(P,X)$ of the Yang--Mills--Higgs functional is gauge equivalent to a smooth pair. Since the estimate (\ref{Lojeq1}) is $\mathcal{G}^{2,2}(P)$ invariant, we may assume in the following that $(A_{\infty},u_{\infty}) \in \mathcal{H}(P,X)$ is smooth. 

The infinitesimal gauge action induces for $s = \pm 1$ the $L^2$-orthogonal splittings
			\begin{align} \label{Lojeq2} H^s(\Sigma, T^*\Sigma \otimes \text{ad}(P)) \oplus H^{s+1}(\Sigma, u_{\infty}^* TX/G) = I^{s+1} \oplus V^{s,{s+1}} \end{align}
with 
	\begin{align*}	I^{s+1} 	&:= \{ (-d_{A_{\infty}} \xi,  L_{u_{\infty}} \xi)\,|\, \xi \in H^{s+1}(\Sigma, \text{ad}(P)) \} \\
									V^{s,s+1} &:= \{ (a,\hat{u}) \in H^s\times H^{s+1}\,|\, d_{A_{\infty}}^* a + L_{u_{\infty}}^* \hat{u} = 0 \}.
	\end{align*}
Define
			$$E: V^{1,2} \rightarrow \mathbb{R}, \qquad E(a,\hat{u}) := \mathcal{YMH}(A_{\infty} + a, \exp_{u_{\infty}} \hat{u}).$$
When $||\hat{u}||_{L^{\infty}}$ is smaller than the injectivity radius of $X$ along $u_{\infty}(P)$, the $L^2$-gradient of $E$ is given by
			$$\nabla E (a,\hat{u}) = \Pi_V \circ T_{(a,\hat{u})} \nabla \mathcal{YMH}(A_{\infty} + a , \exp_{u_{\infty}} \hat{u}) $$
where 
			$$T_{(a,\hat{u})} : T_{(A_{\infty}+a, \exp_{u_{\infty}} \hat{u})} (\mathcal{A}\times \mathcal{S}(P,X)) \rightarrow T_{(A_{\infty},u_{\infty})} (\mathcal{A}\times \mathcal{S}(P,X))$$
												$$T_{(a,\hat{u})}(b, v) := (b, d\exp_{u_0}^{-1} v)$$ 
and $\Pi_V$ denotes the orthogonal projection onto $V^{-1,0}$ in (\ref{Lojeq2}). The next Theorem establishes the {\L}ojasiewicz inequality for $E$ and we show below that this is equivalent to Theorem \ref{ThmLoj1}.

\begin{Theorem} \label{ThmLoj2}
In the setting described above, there exist $\epsilon, C > 0$ and $\gamma \in [\frac{1}{2}, 1)$ such that for all $(a,\hat{u}) \in V^{1,2}$ with $||a||_{H^1} + ||\hat{u}||_{H^2} < \epsilon$ it holds
			\begin{align} \label{Lojeq3} ||\nabla E(a, \hat{u})||_{H^{-1}\times L^2} \geq C | E(a, \hat{u}) - E(0,0) |^{\gamma}. \end{align}
\end{Theorem}

\begin{proof}
$E$ is an analytic functional by assumption \textbf{(A)} and we claim that its Hessian
							$$\nabla^2 E (0,0) : V^{1,2} \rightarrow V^{-1,0}$$
satisfies the elliptic estimate 
					\begin{align} \label{proofLojeq1} ||(b,v)||_{H^1\times H^2} \leq C\left(||\nabla^2 E (0,0) (b,v)||_{H^{-1}\times L^2} + ||(b,v)||_{L^2\times H^1}\right) \end{align}
for all $(b,v)\in V^{1,2}$. We show in Theorem \ref{ThmA} that under these assumptions the {\L}ojasiewicz gradiant inequality is always satisfied and Theorem \ref{ThmLoj2} is a direct consequence of this more general result.

The Hessian $Q$ of the Yang--Mills--Higgs functional at $(A_{\infty},u_{\infty})$ has the shape
					$$Q (b, v) = (d_{A_{\infty}}^*d_{A_{\infty}} b, \nabla_{A_{\infty}}^*\nabla_{A_{\infty}} v) + R(b,v)$$
with some compact operator $R$. As a Hessian, this operator is symmetric and it follows from the gauge-invariance of the Yang--Mills--Higgs functional, that it restricts to the Hessian of $E$, i.e.
					$$\nabla^2 E (0,0) = Q|_{V^{1,2}} : V^{1,2} \rightarrow V^{-1,0}$$
takes indeed values within $V^{-1,0}$. Now consider the operator
			$$\Lambda:= Q + \mathcal{L}_{(A_{\infty},u_{\infty})}\mathcal{L}_{(A_{\infty},u_{\infty})}^*:\, V^{1,2}\oplus I^{2} \rightarrow V^{-1,0}\oplus I^{0}$$
This has the shape			
			$$\Lambda (b,v) = (d_{A_{\infty}}^*d_{A_{\infty}} b + d_{A_{\infty}}d_{A_{\infty}}^* b, \nabla_{A_{\infty}}^*\nabla_{A_{\infty}} v) + \tilde{R}(b,v)$$
with some compact operator $\tilde{R}$. In particular, $\Lambda$ is Fredholm and satisfies the elliptic estimate
					\begin{align} \label{proofLojeq2} ||(b,v)||_{H^1\times H^2} \leq C\left(||\Lambda (b,v)||_{H^{-1}\times L^2} + ||(b,v)||_{L^2\times H^1}\right). \end{align}
Since $\mathcal{L}_{(A_{\infty},u_{\infty})}\mathcal{L}_{(A_{\infty},u_{\infty})}^*$ vanishes on $V^{1,2}$, we have $\Lambda|_{V^{1,2}} = \nabla^2 E$ and (\ref{proofLojeq1}) follows from (\ref{proofLojeq2}).
\end{proof}

\begin{proof}[Proof of Theorem \ref{ThmLoj1}.] \label{proofThmLoj1}
Since the Yang--Mills--Higgs functional is gauge-invariant, it follows from the implicit function theorem that we may assume $(a,\hat{u}) \in V^{1,2}$ with respect to the splitting (\ref{Lojeq2}).
Let $\epsilon > 0$ be sufficiently small, let $(a,\hat{u}) \in V^{1,2}$ with $||a||_{H^1} + ||\hat{u}||_{H^2} < \epsilon$ and denote $(A, u) = (A_{\infty} + a, \exp_{u_{\infty}} \hat{u})$. Then
		\begin{align*}
					||\nabla E(a, \hat{u})||_{H^{-1}\times L^2} &\leq C ||T_{(a,\hat{u})} \nabla \mathcal{YMH} (A,u)||_{H^{-1}\times L^2} \\
																											&\leq C ||\nabla \mathcal{YMH}(A,u)||_{H^{-1}\times L^2}
		\end{align*}
and Theorem \ref{ThmLoj1} follows from Theorem \ref{ThmLoj2}.
\end{proof}

\begin{Remark} Theorem \ref{ThmLoj1} is in fact equivalent to Theorem \ref{ThmLoj2}. To see this denote by $\Pi_I$ the projection onto $I^0$ in (\ref{Lojeq2}). With the same notation as above follows
			\begin{align*}
						&||\Pi_I \circ T_{(a,\hat{u})}\nabla \mathcal{YMH}(A,u) ||_{H^{-1}\times L^2} \\
										&\qquad\leq C ||\mathcal{L}_{(A_{\infty},u_{\infty})}^* T_{(a,\hat{u})} \nabla \mathcal{YMH}(A,u)||_{H^{-2}\times H^{-1}} \\
										&\qquad= C ||(\mathcal{L}_{(A_{\infty},u_{\infty})}^* T_{(a,\hat{u})} - \mathcal{L}_{(A,u)}^*) \nabla \mathcal{YMH}(A,u)||_{H^{-2} \times H^{-1}}
			\end{align*}
where the second equation uses $\mathcal{L}_{(A,u)}^* \nabla \mathcal{YMH}(A,u) = 0$. One verifies that the operator norm of $\mathcal{L}_{(A_{\infty},u_{\infty})}^* T_{(a,\hat{u})} - \mathcal{L}_{(A,u)}^*$ tends to zero as $||a||_{H^1} + ||\hat{u}||_{H^2}$ tends to zero. Hence 
								$$||\Pi_I \circ T_{(a,\hat{u})}\nabla \mathcal{YMH}(A,u) ||_{H^{-1}\times L^2} \leq C(\epsilon) ||\nabla \mathcal{YMH}(A,u)||_{L^2}$$
where $C(\epsilon) > 0$ tends to zero as $\epsilon \rightarrow 0$ and therefore
			\begin{align*} ||\nabla E(a,\hat{u})||_{H^{-1}\times L^2} &\geq  C ||\nabla \mathcal{YMH}(A,u)||_{H^{-1}\times L^2} \end{align*}
for sufficiently small $\epsilon > 0$.
\end{Remark}

\subsection{Interior regularity}

\begin{Theorem} \label{ThmInt}
Let $(A_\infty,u_\infty) \in \mathcal{H}^{1,2}(P,X)$ be a critical point of $\mathcal{YMH}$. There exists $\epsilon_0 > 0$ such that for every $\epsilon \in (0,\epsilon_0)$ there exists $C > 0$ with the following significance: let $T > 1$ and let
				$$(a, \hat{u}): [0,T] \rightarrow H^1(\Sigma, T^*\Sigma\otimes \text{ad}(P))\times H^2(\Sigma, u_{\infty}^* TX/G)$$
be a continuous map such that $(A(t),u(t)) = (A_{\infty} + a(t), \exp_{u_{\infty}} \hat{u}(t)) \subset \mathcal{H}^{1,2}(P,X)$ is a solution of (\ref{floweq1}) and $||a(t)||_{H^1} + ||\hat{u}(t)||_{H^2} < \epsilon$ for all $t \in [0,T]$. Then
				$$\int_{1}^T ||\partial_t a||_{H^1} + ||\partial_t \hat{u}||_{H^2}\, dt \leq C \int_0^T ||\nabla \mathcal{YMH}(A(t),u(t)) ||_{L^2}\, dt.$$
\end{Theorem}

\begin{proof}
By Lemma \ref{Lemma.SobolevReg}, we may assume that $(A_{\infty},u_{\infty}) \in \mathcal{H}(P,X)$ is smooth, after applying a suitable gauge transformation. The idea of the proof is then to show that the derivatives
					$$b := \partial_t a = *d_A(*F_A +  \mu(u))$$
					$$v := \partial_t \hat{u} = d \exp_{u_{\infty}}(\hat{u})^{-1} \left(-\nabla_A^* d_A u - J L_u \mu(u)\right)$$
are solutions of the heat equation up to some perturbation which can be controlled. Note that the later expression is well-defined for sufficiently small $\epsilon > 0$, i.e. when $||\hat{u}||_{L^{\infty}} \leq C ||\hat{u}||_{H^2} \leq C \epsilon$ is smaller than the injectivity radius of $X$ along the image of $u_{\infty}$. Using standard estimates for the heat equation, we can then deduce the interior regularity estimate. In the following fix a small time $0 < t_0 < \min \{1, T/2\}$.

Differentiating $b$ in time gives
		\begin{align*} \partial_t b &= \partial_t \left( *d_A(*F_A +  \mu(u)) \right) = -d_A^* d_A b + [*b, *F_A + \mu(u)] + *d_A (\partial_t \mu(u))
		\end{align*}
From the gauge-invariance follows $d_A^* b + L_u (\partial_t u) = 0$ and hence
		$$\partial_t b + \Delta_A b = [*b, *F_A + \mu(u)] + d_A L_u(\partial_t u) + * d_A(\partial_t \mu(u)).$$
Finally the Bochner-Weizenböck formula yields the relation
			$$\Delta_A b = \nabla_\infty^* \nabla_\infty b + F_A\times b + R_{\Sigma} \times b$$
where $\times$ denotes some bilinear expression and $R_{\Sigma}$ is the Riemann curvature tensor of $\Sigma$. Then follows
		\begin{align*}
				(\partial_t + \nabla_\infty^* \nabla_\infty) b &= \nabla_\infty a \times b + a \times \nabla_\infty b + a\times a \times b + F_\infty \times b + R_{\Sigma} \times b + b\times \mu(u)\\
																					 & \qquad  + a \times L_u (\partial_t u) + a \times (\partial_t \mu(u)) + \nabla_\infty (L_u (\partial_t u)) + \nabla_\infty (\partial_t \mu(u))
		\end{align*}
Now choose a smooth cut-off function $\eta(t)$ such that $\eta(t) = 0$ for $\eta \in [0, t_0/2]$ and $\eta(t) = 1$ for $t \in [t_0, 2t_0]$. Then $\eta b$ satisfies 
			$$(\partial_t + \nabla_\infty^* \nabla_\infty) (\eta b) = \eta (\partial_t + \nabla_\infty^* \nabla_\infty) b + \eta'(t) b.$$
and vanishes at $t=0$. Lemma \ref{LemmaHeat2} shows
		\begin{align*}
				||\eta b||_{L^2([0,2t_0], H^1)} 
									&\leq C t_0^{\frac{1}{4}} ||\eta(\nabla_\infty a \times b + \cdots + b\times \mu(u))||_{L^2([0,2t_0], H^{-\frac{1}{2}})} \\
									& \qquad	+ C t_0^{\frac{1}{4}} ||\eta(a \times L_u (\partial_t u) + \cdots + \nabla_\infty (\partial_t \mu(u)))||_{L^2([0,2t_0], H^{-\frac{1}{2}})} \\
									& \qquad + C t_0^{\frac{1}{2}} ||\eta_t' b||_{L^1([0,2t_0],L^2)}
			\end{align*}
Using the assumption $||a(t)||_{H^1} + ||\hat{u}(t)||_{H^2} < \epsilon$ and the multiplication theorem $H^1\otimes L^2 \rightarrow H^{-\frac{1}{2}}$ it follows:
		$$||\eta(\nabla_\infty a \times b + \cdots + b\times \mu(u))||_{L^2([0,2t_0], H^{-\frac{1}{2}})} \leq C\epsilon ||\eta b||_{L^2([0,2t_0], H^1)}$$
		$$||\eta(a \times L_u (\partial_t u) + \cdots + \nabla_\infty (\partial_t \mu(u)))||_{L^2([0,2t_0], H^{-\frac{1}{2}})} \leq C \epsilon ||\eta v||_{L^2([0,2t_0], H^2)}$$
By choosing $t_0$ sufficiently small we thus obtain
		\begin{align} \label{IReq1}
			||\eta b||_{L^2([0,2t_0], H^1)} \leq  c ||\eta' b||_{L^1([0,2t_0], L^2)} + C t_0^{\frac{1}{4}} ||\eta v||_{L^2([0,2t_0], H^2)}
		\end{align}	

Next, we need to obtain a similar estimate for $v$. Define
		$$\Psi : H^1(\Sigma, T^*\Sigma\otimes \text{ad}(P))\times H^2(\Sigma, u_{\infty}^*TX/G) \rightarrow H^{-\frac{1}{2}}(\Sigma, u_{\infty}^* TX/G)$$
		$$\Psi(a,\hat{u}) := d\exp_{u_{\infty}}^{-1} ( - \nabla_A^* d_A u - J L_u \mu(u))$$
with $A = A_{\infty} + a$ and $u = \exp_{u_{\infty}} \hat{u}$.  This is continuously differentiable and satisfies $\Psi(0,0) = 0$. In particular, 
		\begin{align*}
				v &= \Psi(a, \hat{u}) = d\Psi(0,0) (a, \hat{u}) + q(a, \hat{u})
		\end{align*}
where $q$ vanishes to the first order. Differentiating this equation with respect to $t$ yields
		\begin{align*} \partial_t v &= d\Psi(0,0) (b, v) + \partial_t q(a,\hat{u}) \\
																&= -\nabla_{A_{\infty}}^* \nabla_{A_{\infty}} v - \nabla_{v} (J L_{u_{\infty}} \mu(u_{\infty})) - \nabla_{A_{\infty}}^* L_{u_{\infty}} b + \nabla b \times d_{A_\infty} u_{\infty} \\
																&\qquad + d q(a_t, \hat{u}_t) [b, v]
		\end{align*}
Let $\eta$ be the same cut-off function as above. Then $\eta v$ vanishes at $t=0$ and solves the equation
		\begin{align*}
				(\partial_t + \nabla_{A_{\infty}}^* \nabla_{A_{\infty}}) (\eta v) &= - \nabla_{\eta v} (J L_{u_{\infty}} \mu(u_{\infty}) - \nabla_{A_{\infty}}^* L_{u_{\infty}} (\eta b) + \nabla (\eta b) \times d_{\infty} u_{\infty} \\
									& \qquad + \eta \partial_t q(a_t, \hat{u}_t) [b, v] + \eta' v
		\end{align*}
It follows from Lemma \ref{LemmaHeat2} that
		\begin{align*}
				||\eta v||_{L^2([0,2t_0], H^1)} 
									&\leq C t_0^{\frac{1}{4}} ||- \nabla_{\eta v} (J L_{u_{\infty}} \mu(u_{\infty}) - \ldots + \nabla (\eta b) \times d_{\infty} u_{\infty}||_{L^2([0,2t_0], H^{-\frac{1}{2}})} \\
									& \qquad	+ C t_0^{\frac{1}{4}} ||\eta d q(a_t, \hat{u}_t) [b, v]||_{L^2([0,2t_0], H^{-\frac{1}{2}})} \\
									& \qquad + C t_0^{\frac{1}{2}}||\eta_t' v||_{L^1([0,2t_0],L^2)}
			\end{align*}
The first term satisfies the estimate
		\begin{align*}
				&||- \nabla_{\eta v} (J L_{u_{\infty}} \mu(u_{\infty}) - \cdots + \nabla (\eta b) \times d_{\infty} u_{\infty}||_{L^2([0,2t_0], H^{-\frac{1}{2}})} \\
						&\qquad\leq C ||\eta v||_{L^2([0,2t_0], H^1)} + C ||\eta b||_{L^2([0,2t_0], H^1)}
		\end{align*}				
and it follows from the definition of $q$ that
			$$||\eta d q(a_t, \hat{u}_t) [b, v]||_{L^2([0,2t_0], H^{-\frac{1}{2}})} \leq C \epsilon \left( ||\eta b||_{L^2([0,2t_0], H^1)} + ||\eta v||_{L^2([0,2t_0], H^2)} \right)$$
For sufficiently small $t_0 > 0$ we thus get
		\begin{equation} \label{IReq2}
		\begin{aligned} 
			||\eta v||_{L^2([0,2t_0], H^1)} &\leq  C ||\eta'_t v||_{L^1([0,2t_0], L^2)} \\
																			&\qquad + C t_0^{\frac{1}{4}} (||\eta v||_{L^2([0,2t_0], H^2)} + ||\eta b||_{L^2([0,2t_0], H^1)})
		\end{aligned}
		\end{equation}	

Finally, differentiating the holomorphicity condition $\bar{\partial}_A u = \bar{\partial} u + L_u A^{0,1} = 0$, we obtain by elliptic regularity the estimate
				\begin{align} \label{IReq3} ||v_t||_{H^2} \leq C (||v_t||_{H^1} + ||b_t||_{H^1}). \end{align}
Combining (\ref{IReq1}, (\ref{IReq2}) and (\ref{IReq3}) yields for sufficiently small $t_0 > 0$
		\begin{align*}
			||\eta b||_{L^2([0,2t_0], H^1)} + ||\eta v||_{L^2([0,2t_0], H^2)} \leq  C (||\eta' b||_{L^1([0,2t_0], L^2)} + ||\eta' v||_{L^1([0,2t_0], L^2)} ).
		\end{align*}		
In particular
		\begin{align*}
			||b||_{L^1([t_0,2t_0], H^1)} + ||\eta v||_{L^1([t_0,2t_0], H^2)} 
							&\leq   t_0^{\frac{1}{2}} (||\eta b||_{L^2([0,2t_0], H^1)} + ||\eta v||_{L^2([0,2t_0], H^2)}) \\
							&\leq C t_0^{\frac{1}{2}} (||\eta' b||_{L^1([0,2t_0], L^2)} + ||\eta' v||_{L^1([0,2t_0], L^2)} ) \\
							&\leq C t_0^{-\frac{1}{2}} (||b||_{L^1([0,2t_0], L^2)} + ||v||_{L^1([0,2t_0], L^2)} ).
		\end{align*}
The proof follows now by subdividing the interval $[0,T]$ into smaller intervals of length $t_0$ and applying the estimate above to each pair of successive subintervals.

\end{proof}

\subsection{The convergence theorem}

Theorem \ref{ThmConv0} is a slightly weaker version of the next theorem.

\begin{Theorem} \label{ThmConv}
Assume \textbf{(A)}, \textbf{(B)} and \textbf{(C)}. Let $(A_0, u_0) \in \mathcal{H}^{1,2}(P,X)$ and let $(A(t), u(t))$ be the solution of (\ref{floweq2}). There exist $C, \beta > 0$ such that for all $T > 0$
							$$\int_T^{\infty} ||\partial_t A(t)||_{H^1} + ||\partial_t u(t)||_{H^2}\, dt \leq C T^{-\beta}.$$
In particular, $(A(t),u(t))$ converges uniformly in $H^1\times H^2$ to a critical point $(A_{\infty}, u_{\infty})$ of $\mathcal{YMH}$.
\end{Theorem}

The following compactness result arises from a combination of Gromov compactness for holomorphic curves and Uhlenbeck compactness.

\begin{Proposition} \label{PropCpct}
Assume \textbf{(B)}, \textbf{(C)} and let $(A(t), u(t)) \subset \mathcal{H}^{1,2}(P,X)$ be a solution of (\ref{floweq2}). Then there exists sequences of times $t_j \rightarrow \infty$ and of gauge transformations $k_j \in \mathcal{G}^{2,2}(P)$ and a critical point $(A_{\infty}, u_{\infty})$ of $\mathcal{YMH}$ such that $k_j(A(t_j), u(t_j))$ converges to $(A_{\infty}, u_{\infty})$ in $H^1 \times H^2$.
\end{Proposition}

\begin{proof}
This is a special case of Theorem 1.2 in \cite{Venugopalan:2016}. 
\end{proof}

\begin{Remark}
In general, we expect the convergence of $u(t_j)$ only modulo bubbling in finitely many fibers, as stated in \cite{Venugopalan:2016} Theorem 1.2.  Assumption \textbf{(B)} rules the formation of bubbles out and is crucial for the result stated in Proposition \ref{PropCpct}.
\end{Remark}

\begin{proof}[Proof of Theorem \ref{ThmConv}]
Let $(A(t), u(t))\subset \mathcal{H}^{1,2}(P,X)$ be a solution of (\ref{floweq2}). Let $t_j \rightarrow \infty$, $k_j \in \mathcal{G}^{2,2}(P)$ and $(A_{\infty}, u_{\infty})$ be as in Proposition \ref{PropCpct} above. Choose $\epsilon > 0$, such that the {\L}ojasiewicz gradient inequality in Theorem \ref{ThmLoj1} is satisfied with respect to $(A_{\infty}, u_{\infty})$. Let $\delta \in (0, \epsilon)$ and choose $j \geq 1$ such that
							$$||A_{\infty} - k_j A(t_j)||_{H^1} + ||\exp_{u_{\infty}}^{-1} (k_j u(t_j))||_{H^2} < \delta.$$
Since the gradient flow and the {\L}ojasiewicz inequality are $\mathcal{G}^{2,2}(P)$-equivariant, we may assume $k_j = \mathds{1}$ and $t_j = 0$.

The gradient flow depends continuously in the $C^0(H^1\times H^2)$ topology on the initial conditions by Theorem \ref{ThmVen}. Since the flow is constant at the critical point $(A_{\infty}, u_{\infty})$ this yields
				$$||A(1) - A(0)||_{H^1} + ||\exp_{u_{\infty}}^{-1} u_1 - \exp_{u_{\infty}}^{-1} u_0||_{H^2} \leq \rho(\delta)$$
where $\rho(\delta) \rightarrow 0$	as $\delta \rightarrow 0$. Define 
			$$\overline{T} := \inf \{ t > 0 \,|\, ||A(t) - A_{\infty}||_{H^1} + ||\exp_{u_{\infty}}^{-1} u_t||_{H^2} \geq \epsilon \}.$$
By choosing $\delta > 0$ sufficiently small, we can guarantee $\overline{T} > 1$. For $1 < s < \overline{T}$ define $\hat{u}(s) := \exp_{u_{\infty}}^{-1} u(s)$. The interior regularity estimate in Theorem \ref{ThmInt} and the {\L}ojasiewicz gradient inequality in Theorem \ref{ThmLoj1} yield
		\begin{align*}
				||A(s) - A_{\infty}||_{H_1} + ||\hat{u}(s) ||_{H^2}
						&\leq \rho(\delta) + \int_1^s ||\partial_t A(t)||_{H^1} + ||\partial_t \hat{u}(t)||_{H^2}\, dt \\
						&\leq \rho(\delta) + C \int_0^s ||\partial_t A(t)||_{L^2} + ||\partial_t u(t)||_{L^2}\, dt \\
						&\leq \rho(\delta) + C \int_0^s \frac{||\nabla \mathcal{YMH}(A,u)||^2_{L^2\times L^2}}{(\mathcal{YMH}(A,u) - \mathcal{YMH}(A_{\infty},u_{\infty}))^{\gamma}}\, dt \\
						&\leq	\rho(\delta) + C\left(\mathcal{YMH}(A(0),u(0)) - \mathcal{YMH}(A_{\infty}, u_{\infty})\right)^{1-\gamma}.
		\end{align*}
For $\delta > 0$ sufficiently small, this shows $\overline{T} = \infty$ and the integral $\int_1^\infty ||\partial_t A(t)||_{H^1} + ||\partial_t u(t)||_{H^2} \,dt < \infty$ is finite.
This proves that $(A(t), u(t))$ converges uniformly in $H^1\times H^2$ to a critical point $(\tilde{A}_{\infty}, \tilde{u}_{\infty})$ of the Yang--Mills--Higgs functional. 

Repeating the argument from above, with respect to the critical point $(\tilde{A}_{\infty},\tilde{u}_{\infty})$ we obtain for all sufficently large $T$
		$$\int_T^\infty ||\partial_t A(t)||_{H^1} + ||\partial_t \hat{u}(t)||_{H^2} \,dt  \leq f(T - 1)^{1-\gamma} $$
with $f(t) := (\mathcal{YMH}(A(t), u(t)) - \mathcal{YMH}(\tilde{A}_{\infty}, \tilde{u}_{\infty}))$. Since
		\begin{align*}
				f'(t) = - ||\nabla \mathcal{YMH}(A(t),u(t))||_{L^2}^2  \leq - C f(t)^{2\gamma} 
		\end{align*}
it follows $f(t) \leq C t^{\frac{1}{1 - 2\gamma}}$ and hence
				$$\int_T^\infty ||\partial_t A(t)||_{H^1} + ||\partial_t \hat{u}(t)||_{H^2} \,dt \leq C (T-1)^{\frac{1 - \gamma}{1- 2\gamma}}$$
for all sufficiently large $T$. This is equivalent to the estimate in the Theorem and completes the proof.
\end{proof}

We state some consequences of the proof for later reference.

\begin{Cor} \label{CorLoj}
Assume \textbf{(A)}, \textbf{(B)}, \textbf{(C)} and let $(B,v) \in \mathcal{H}(P,X)$ be a critical point of the Yang--Mills--Higgs functional. There exist $C, \epsilon_0 > 0$ and $\gamma \in [\frac{1}{2},1)$ with the following significance: let $(A, u) : [0,\infty) \rightarrow \mathcal{H}(P,X)$ be a solution of (\ref{floweq2}) satisfying $||A(0) - B||_{H^1} + ||\exp_{v}^{-1} u(0)||_{H^2} < \epsilon_0$ and $\mathcal{YMH}(A(t),u(t)) \geq \mathcal{YMH}(B,v)$ for all $t > 0$. Then
				\begin{enumerate}
						\item The limit satisfies $\mathcal{YMH}(A_{\infty}, u_{\infty}) = \mathcal{YMH}(B, v)$.
						\item For every $\epsilon > 0$ exists $\delta \in (0,\epsilon_0)$ such that 
												$$\int_0^{\infty}||\partial_t A(t)||_{H^1} + ||\partial_t u(t)||_{H^2} \, dt < \epsilon$$
									if $||A(0) - B||_{H^1} + ||\exp_{v}^{-1} u(0)||_{H^2} < \delta$.
				\end{enumerate}
\end{Cor}

\section{Uniqueness and the Kempf--Ness theorem}

\subsection{Uniqueness of critical points}

The next result is a reformulation of Theorem \ref{ThmB} in the introduction and the analogue of the Ness uniqueness theorem in finite dimensional GIT. The proof is based on arguments of Chen--Sun \cite{ChenSun:2010} in the finite dimensional differentiable setting.

\begin{Theorem}[\textbf{Uniqueness of critical points}] \label{ThmNU}
Assume \textbf{(A)}, \textbf{(B)} and \textbf{(C)}. Let $(A_0,u_0) \in \mathcal{H}^{1,2}(P,X)$ and $(A_{\infty}, u_{\infty})$ be the limit of the Yang--Mills--Higgs flow (\ref{floweq2}) starting at $(A_0,u_0)$. Then $(A_{\infty}, u_{\infty}) \in \overline{(\mathcal{G}^c)^{2,2} (A_0, u_0)}$ (the $H^1\times H^2$ closure) and 
					$$\mathcal{YMH}(A_{\infty}, u_{\infty}) = \inf_{g \in (\mathcal{G}^c)^{2,2}(P)} \mathcal{YMH}(gA_0, g u_0).$$
Moreover, if $(B,v) \in \overline{(\mathcal{G}^c)^{2,2} (A_0,u_0)}$ and $\mathcal{YMH}(B, v) = \mathcal{YMH}(A_{\infty}, u_{\infty})$, then $(B, v) \in \mathcal{G}^{2,2}(A_{\infty}, u_{\infty})$.				
\end{Theorem}

\begin{proof}
The proof consists of four steps.\\

\textbf{Step 1:} \textit{Let $(B,v) \in \mathcal{H}^{1,2}(P,X)$ and let $g_j : [0,\infty) \rightarrow (\mathcal{G}^c)^{2,2}(P)$ satisfy
															$$g_j(t)^{-1} \dot{g}_j(t) = \textbf{i} \left(F_{g_j(t)^{-1}B} + \mu(g_j(t)^{-1} v)\right)$$
								for $j \in \{0,1\}$. Using the Cartan decomposition, write
											$$g_1(t) = g_0(t) e^{\textbf{i} \eta(t)} k(t)$$
								with $\eta(t) \in H^2(\Sigma, \text{ad}(P))$ and $k(t) \in \mathcal{G}^{2,2}(P)$. Then $\eta(t)$ and $k(t)$ are uniformly bounded in $H^2$.
								}\\
								
Denote by $\pi : G^c \rightarrow G^c/G$ the canonical projection. The homogeneous space $G^c/G$ is a complete Riemannian manifold with nonpositive curvature and for $t>0$ the curve $\gamma(s,t) := \pi(g_0(t) e^{\textbf{i}s \eta(t)})$ is pointwise the unique geodesic of length $||\eta(t)||$ connecting $\pi(g_0(t))$ and $\pi(g_1(t))$. This yields
		\begin{align*}
				\partial_t ||\eta(t)||^2 &= 2 \int_0^1 \langle \nabla_t \partial_s \gamma, \partial_s \gamma \rangle\, ds 
																			= 2 \int_0^1 \partial_s \langle \partial_t \gamma, \partial_s \gamma \rangle \, ds \\
																 &= 2 \langle g_1(t)^{-1}\dot{g}_1(t), \textbf{i}\eta(t) \rangle - \langle g_0(t)^{-1}\dot{g}_0(t), \textbf{i} \eta(t) \rangle \\
															   &= 2\langle *F_{g_1(t)^{-1} B} - *F_{g_0(t)^{-1} B}, \eta(t) \rangle + 2 \langle \mu(g_1(t)^{-1}v) - \mu(g_0(t)^{-1} v) , \eta(t) \rangle.
		\end{align*}
Abbreviate $(B_{s,t}, v_{s,t}) := e^{-\textbf{i}s\eta(t)} g_0(t)^{-1}(B, v)$. Then
		\begin{align*}
				\partial_t ||\eta(t)||^2 &= 2 \int_0^1	\partial_s \langle *F_{B_{s,t}} + \mu(v_{s,t}), \eta(t) \rangle \, ds \\
																 &= -2  \langle \Delta_{B_{s,t}} \eta + L_{v_{s,t}}^*L_{v_{s,t}} \eta(t) , \eta(t) \rangle \\
																 &= - \Delta (||\eta(t)||^2 ) - 2 \int_0^1 \left( ||L_{u_{s,t}} \eta(t)||^2 + ||d_{A_{s,t}} \eta(t)||^2 \right)\, ds
		\end{align*}
Thus $||\eta||^2$ satisfies the differential inequality $(\partial_t + \Delta) ||\eta||^2 \leq 0$ and by the maximum principle for the heat equation $\eta(t)$ is uniformly bounded in $L^{\infty}$. Since $(B_j(t),v_j(t)) := (g_j(t)^{-1} B, g_j(t)^{-1} v)$ satisfies (\ref{floweq2}), it converge uniformly in $H^1\times H^2$ by Theorem \ref{ThmConv}. Hence it follows from the equation
					$$ B_1(t)= g_1(t)^{-1} B = k(t)^{-1}e^{-\textbf{i} \eta(t)} g_0(t)^{-1} B = k(t)^{-1}e^{-\textbf{i} \eta(t)} B_0(t)$$
and elliptic bootstrapping that $\eta(t)$ and $k(t)$ are uniformly bounded in $H^2$.\\

\textbf{Step 2:} \textit{Let $(B_0, v_0), (B_1,v_1) \in \mathcal{H}^{1,2}(P,X)$ be critical points of the Yang--Mills--Higgs functional. If $(B_1,v_1) \in (\mathcal{G}^c)^{2,2}(B_0,v_0)$, then $(B_1,v_1) \in \mathcal{G}^{2,2}(B_0, v_0)$.}\\

Choose $\tilde{g} \in (\mathcal{G}^c)^{2,2}(P)$ such that $\tilde{g}^{-1}(B_1,v_1) = (B_0, v_0)$. By Theorem \ref{ThmVen} there exist $g_0, g_1: [0,\infty) \rightarrow (\mathcal{G}^c)^{2,2}(P)$ solving
					$$g_0^{-1} \dot{g}_0 = *F_{B_0} + \mu(v_0), \qquad g_0(0) = \mathds{1}, \qquad g_0^{-1}(t)(B_0,v_0) = (B_0, v_0)$$
					$$g_1^{-1} \dot{g}_1 = *F_{B_1} + \mu(v_1), \qquad g_1(0) = \tilde{g}, \qquad  g_1^{-1}(t)(B_0,v_0) = (B_1,v_1).$$
Both of these curves satisfy the conditions of Step 1 with $(B,v) = (B_0,v_0)$. Using the same notation as in Step 1, write $g_1(t) = g_0(t) e^{\textbf{i} \eta(t)} k(t)$, and conclude that there exists a sequence $t_j \rightarrow \infty$ such that $k(t_j) \rightarrow k_{\infty}$ and $\eta(t_j) \rightarrow \eta_{\infty}$ converge weakly in $H^2$ and strongly in $W^{1,p}$. Then
			$$B_{s,t_j} \stackrel{L^p}{\longrightarrow} B_{s,\infty} :=e^{-\textbf{i}s \eta_{\infty}} B_0, \qquad d_{B_{s,t_j}} \eta(t_j) \stackrel{L^p}{\longrightarrow} d_{B_{s,\infty}} \eta_{\infty}$$
where $(B_{s,t}, v_{s,t}) := e^{-\textbf{i}s\eta(t)} (B_0,v_0)$ as in Step 1. It follows from the calculation in Step 1 that
				$$\partial_t ||\eta||_{L^2}^2 = -2 \int_0^1 ||d_{B_{s,t}} \eta||_{L^2}^2 + ||L_{v_{s,t}} \eta||_{L^2}^2\, ds$$ 
and we may assume in addition
				$$\lim_{j\rightarrow \infty} ||d_{B_{s,t_j}} \eta(t_j)||_{L^2} + ||L_{v_{s,t_j}} \eta(t_j)||_{L^2}  = 0.$$
At $s= 0$ we obtain $\mathcal{L}_{(B_0,v_0)} \eta_{\infty} = 0$ and hence
			$$(B_1,v_1) = g_1(t_j)^{-1}(B_1,v_1) \stackrel{L^p}{\longrightarrow}  k_{\infty}^{-1} e^{-\textbf{i} \eta_{\infty}} (B_0, v_0) = k_{\infty}^{-1} (B_0, v_0).$$
This shows $(B_1,v_1) \in \mathcal{G}^{2,2}(B_0, v_0)$ and completes the proof of Step 2. \\

\textbf{Step 3:} \textit{Let $(A_0,u_0), (B_0,v_0) \in \mathcal{H}^{1,2}(P,X)$ and denote by $(A_{\infty},u_{\infty}), (B_{\infty}, v_{\infty})$ the limits of the Yang--Mills--Higgs flow (\ref{floweq2}) starting at $(A_0,u_0)$, $(B_0,v_0)$ respectively. If $(B_0,v_0) \in (\mathcal{G}^c)^{2,2}(A_0,u_0)$, then $(B_{\infty},v_{\infty}) \in \mathcal{G}^{2,2}(A_{\infty}, u_{\infty})$.}\\

Denote by $(A(t),u(t))$ and $(B(t),v(t))$ the solutions of (\ref{floweq2}) starting at $(A_0,u_0)$ and $(B_0,v_0)$ respectively. Choose $\tilde{g} \in (\mathcal{G}^c)^{2,2}(P)$ such that $(B_0,v_0) = \tilde{g}^{-1} (A_0,u_0)$. Then, by Theorem \ref{ThmVen}, there exist $g_0, g_1: [0,\infty) \rightarrow (\mathcal{G}^c)^{2,2}(P)$ solving
			$$g_0^{-1}\dot{g_0}(t) = \Phi(A(t), u(t)), \qquad g_0(0) = \mathds{1}, \qquad g_0^{-1}(t)(A_0,u_0) = (A(t), u(t))$$
			$$g_1^{-1}\dot{g_1}(t) = \Phi(B(t),u(t)), \qquad g_1(0) = \tilde{g}, \qquad g_1^{-1}(t)(A_0,u_0) = (B(t), v(t)).$$
Both of these curves satisfy the conditions of Step 1 with $(B,v) = (A_0,u_0)$. Using the same notation as in Step 1, write $g_1(t) = g_0(t) e^{\textbf{i}\eta(t)}k(t)$. Then there exists a sequence $t_j \rightarrow \infty$ such that $\eta(t_j) \rightarrow \eta_{\infty}$ and $k(t_j) \rightarrow k_{\infty}$ converge weakly in $H^2$ and strongly in $W^{1,p}$. As $j$ tends to infinity in the equation
					$$(B(t_j),v(t_j)) = k(t_j)^{-1} e^{-\textbf{i}\eta(t_j)} (A(t_j),u(t_j))$$
both sides converge in $L^p\times W^{1,p}$ and this yields $(B_{\infty}, v_{\infty})=  k_{\infty}^{-1} e^{-\textbf{i}\eta_{\infty}} (A_{\infty},u_{\infty})$. This proves $(\mathcal{G}^c)^{2,2}(A_{\infty}, u_{\infty}) = (\mathcal{G}^c)^{2,2}(B_{\infty},v_{\infty})$ and Step 3 follows from Step 2. \\

\textbf{Step 4:} \textit{If $(B,v) \in \overline{(\mathcal{G}^c)^{2,2} (A_0,u_0)}$ and $\mathcal{YMH}(B, v) = \mathcal{YMH}(A_{\infty}, u_{\infty})$, then $(B,v) \in \mathcal{G}^{2,2}(A_{\infty}, u_{\infty})$.}\\

It follows from Step 3 that
			$$\mathcal{YMH}(A_{\infty}, u_{\infty}) = \inf_{g\in (\mathcal{G}^c)^{2,2}(P)} \mathcal{YMH}(g A_0, g u_0) =: m.$$
Note that the solution $(B(t),v(t))$ of the Yang--Mills--Higgs flow (\ref{floweq2}) starting at $(B,v)$ remains in the closure $\overline{(\mathcal{G}^c)^{2,2} (A_0,u_0)}$. Hence $\mathcal{YMH}(B(t),v(t)) = m$ is constant, $(B(t),v(t))$ a constant flow line and $(B,v)$ a critical point.

Choose $(A^{(j)},u^{(j)}) \in (\mathcal{G}^c)^{2,2} (A_0,u_0)$ converging in $H^1\times H^2$ to $(B,v)$ and denote the limit of the Yang--Mills--Higgs flow starting at $(A^{(j)},u^{(j)})$ by $(B^{(j)},v^{(j)})$. Corollary \ref{CorLoj} shows that $(B^{(j)},v^{(j)})$ converges to $(B,v)$ in the $H^1\times H^2$ topology. By Step 3, there exists $k_j \in \mathcal{G}^{2,2}(P)$ such that $(B^{(j)}, v^{(j)}) = (k_j A_{\infty}, k_j u_{\infty})$. Since the connections $B^{(j)}$ are uniformly bounded in $H^1$, the gauge transformations $k_j$ are uniformly bounded in $H^2$ and after passing to a subsequence, we may assume that $k_j \rightarrow k_{\infty}$ converges weakly in $H^2$ and strongly in $W^{1,p}$. It follows $(B, v) = (k_{\infty} A_{\infty}, k_{\infty} u_{\infty})$ and this completes the proof.

\end{proof}

\begin{Theorem} \label{ThmLimitStability}
Assume \textbf{(A)}, \textbf{(B)} and \textbf{(C)}. Let $(A,u) \in \mathcal{H}^{1,2}(P,X)$ and denote by $(A_{\infty}, u_{\infty})$ the limit of the flow (\ref{floweq0}) starting at $(A,u)$. Then
			\begin{enumerate}
					\item $(A,u)$ is stable if and only if $\mathcal{L}_{(A_{\infty},u_{\infty})}$ is injective.
					\item $(A,u)$ is polystable if and only if $(A_{\infty},u_{\infty}) \in \mathcal{G}^c(A,u)$ and $*F_{A_{\infty}} + \mu(u_{\infty}) = 0$.
					\item $(A,u)$ is semistable if and only if $*F_{A_{\infty}} + \mu(u_{\infty}) = 0$.
					\item $(A,u)$ is unstable if and only if $*F_{A_{\infty}} + \mu(u_{\infty}) \neq 0$.
			\end{enumerate}
\end{Theorem}

\begin{Remark}
We call an element $(A,u) \in \mathcal{H}^{1,2}(P,X)$ stable, polystable, semistable or unstable, if every smooth element of $(\mathcal{G}^c)^{2,2}(A,u)$ is stable in the sense of Definition \ref{DefStability}. Note that for every stable pair $(A,u)$ the extension of the infinitesimal action
				$$\mathcal{L}_{(A,u)}^c : H^2(\Sigma,\text{ad}(P)) \rightarrow H^1(\Sigma,T^*\Sigma \otimes \text{ad}(P)) \oplus H^2(\Sigma, u^* TX/G)$$
remains injective. This follows from Lemma \ref{Lemma.SobolevReg} and elliptic regularity.
\end{Remark}

\begin{proof}
The unstable, semistable and polystabe characterization follow directly from Theorem \ref{ThmNU}.

For the stable case, note that every stable orbit has discrete $\mathcal{G}^c$-isotropy and this proves one direction. Conversely, the limit satisfies the critical point equation $\mathcal{L}_{(A_{\infty},u_{\infty})}(*F_{A_{\infty}} + \mu(u_{\infty})) = 0$. Hence, when $\mathcal{L}_{(A_{\infty},u_{\infty})}$ is injective, $(A_{\infty},u_{\infty})$ is stable. Since the subset of stable pairs $\mathcal{H}_s^{1,2}(P,X)$ is open by Proposition \ref{PropOpen} below, this implies $(A(t),u(t)) \in \mathcal{H}_s^{1,2}(P,X)$ for sufficiently large $t$ and $(A,u)$ is stable.
\end{proof}

\begin{Proposition} \label{PropOpen}
Assume \textbf{(A)}, \textbf{(B)} and \textbf{(C)}. The subsets of stable and semistable pairs $\mathcal{H}_{ss}^{1,2}(P,X)  \subset \mathcal{H}_s^{1,2}(P,X) \subset \mathcal{H}^{1,2}(P,X)$ are open subsets in the $H^1\times H^2$-topology.
\end{Proposition}

\begin{proof}
The semistable case follows from Corollary \ref{CorLoj}. The stable case follows from a suitable application of the implicit function theorem: Suppose $(A,u) \in \mathcal{H}^{1,2}(P,X)$ solves the vortex equation $*F_A + \mu(u) = 0$ and $\mathcal{L}_{(A,u)}$ is injective. Then $\mathcal{L}^c_{(A,u)}$ is also injective and, since
				$$\langle \mathcal{L}_{(A,u)}^c \textbf{i}\xi , (a,\hat{u}) \rangle_{L^2\times L^2} = \langle \xi, *d_A a + d \mu(u) \hat{u} \rangle_{L^2}$$
for all $\xi \in H^2(\Sigma, \text{ad}(P))$ and $(a,\hat{u}) \in T_{(A,u)} (\mathcal{A}(P)\times \mathcal{S}(P,X))$, $(A,u)$ is a regular point for the moment map $\Phi(A,u) = *F_A + \mu(u)$. It follows that 
					$$\mathcal{Z} := \Phi^{-1}(0) \subset \mathcal{A}^{1,2}(P)\times \mathcal{S}^{2,2}(P,X)$$
is a submanifold locally around $(A,u)$ and the orthogonal complement of $T_{(A,u)} \mathcal{Z}$ coincides with the image of $H^2(\Sigma, \textbf{i} \,\text{ad}(P))$ under $\mathcal{L}_{(A,u)}^c$. Hence
				$$H^2(\Sigma, T^*\Sigma \otimes \text{ad}(P^c))\times \mathcal{Z} \rightarrow \mathcal{A}^{1,2}(P)\times\mathcal{S}^{2,2}(P,X), \qquad (\xi, z) \mapsto \exp(\textbf{i} \xi)z$$
restricts to a diffeomorphism between neighborhoods of $(0;(A,u))$ and $(A,u)$. In particular, $(A,u)$ is an interior point of $\mathcal{H}_s^{1,2}(P,X)$.

\end{proof}

\subsection{Kempf--Ness theorem}

Let $(A,u) \in \mathcal{H}^{1,2}(P,X)$ and define the $1$-form $\alpha_{(A,u)}: T(\mathcal{G}^c)^{2,2}(P) \rightarrow \mathbb{R}$ by
	\begin{align} \label{KNeq1} \alpha_{(A,u)}(g,\hat{g}) := - \int_{\Sigma} \left\langle *F_{g^{-1}A} + \mu(g^{-1} u); \text{Im}(g^{-1} \hat{g}) \right\rangle \, dvol_{\Sigma}. \end{align}
It is straight forward to check that $\alpha_{(A,u)}$ is exact, $\mathcal{G}^{2,2}(P)$-invariant and integrates to a unique $\mathcal{G}^{2,2}(P)$-invariant functional
	\begin{align} \label{KNeq2}\Psi_{(A,u)} : (\mathcal{G}^c)^{2,2}(P) \rightarrow \mathbb{R} \end{align}
satisfying $\Psi_{(A,u)}(\mathds{1}) = 0$ (see e.g. \cite{Mundet:2000}). We call $\Psi_{(A,u)}$ the Kempf--Ness functional associated to $(A,u)$. 

Theorem \ref{ThmVen} shows that for every $g_0 \in (\mathcal{G}^c)^{2,2}(P)$ the negative gradient flow
	\begin{align} \label{KNeq3} g^{-1}(t)\dot{g}(t) = - g^{-1}(t)\nabla \Psi_{(A,u)}(g(t)) = - \textbf{i} (*F_{g(t)^{-1} A} + \mu(g(t)^{-1}u)) \end{align}
has a unique solution $g \in C^0([0,\infty), (\mathcal{G}^c)^{2,2}(P))$ satisfying $g(0) = g_0$. This flow intertwines with the Yang--Mills--Higgs flow in the following sense
			$$\text{$g(t)$ solves (\ref{KNeq3}) } \quad \Longrightarrow \quad \text{$(A(t),u(t)) := (g(t)^{-1}A, g(t)^{-1} u)$ solves (\ref{floweq2})}.$$
We will repetitively make use of the fact that $\Psi_{(A,u)}$ is convex along geodesics in $(\mathcal{G}^c)^{2,2}(P)/\mathcal{G}^{2,2}(P)$. This amounts to the formula
	\begin{align} \label{KNeq4} \frac{d^2}{dt^2} \Psi_{(A,u)}( g e^{\textbf{i}t\xi}) = \left|\left| \mathcal{L}_{e^{-\textbf{i}t\xi}g^{-1}(A,u)} \xi \right|\right|_{L^2}^2 \geq 0\end{align}
for $g \in (\mathcal{G}^c)^{2,2}(P)$ and $\xi \in H^2(\Sigma, \text{ad}(P))$.

The Kempf--Ness theorem relates the stability of the pair $(A,u)$ to global properties of the functional $\Psi_{(A,u)}$. The stable case is due to Mundet \cite{Mundet:2000}, see Remark \ref{RmkStableKN} below. The remaining cases are the content of the next theorem which is a reformulation of Theorem \ref{ThmC} in the introduction.

\begin{Theorem} \label{ThmKN}
Assume \textbf{(A)}, \textbf{(B)} and \textbf{(C)} and let $(A,u) \in \mathcal{H}^{1,2}(P,X)$.
		\begin{enumerate}
				\item $(A,u)$ is polystable if and only if $\Psi_{(A,u)}$ has a critical point.
				\item $(A,u)$ is semistable if and only if $\Psi_{(A,u)}$ is bounded below.
				\item $(A,u)$ is unstable if and only if $\Psi_{(A,u)}$ is unbounded below.								
		\end{enumerate}
\end{Theorem}

\begin{proof}
The polystable case follows from (\ref{KNeq1}). For the other two cases let $g_0 \in (\mathcal{G}^c)^{2,2}(P)$ and $g: [0,\infty) \rightarrow (\mathcal{G}^c)^{2,2}(P)$ be the solution of (\ref{KNeq3}) starting at $g_0$. Then
			$$\frac{d}{dt} \Psi_{(A,u)}(g(t)) = \alpha_{(A,u)}(g(t), \dot{g}(t)) = - ||*F_{g(t)^{-1}A} + \mu(g(t)^{-1}u) ||_{L^2}^2.$$
If $(A,u)$ is unstable, Theorem \ref{ThmNU} shows that the right hand side is bounded above by a strictly negative constant and hence $\Psi_{(A,u)}$ is unbounded below. Conversely, assume that $(A,u)$ is semistable. Then $(A(t),u(t)) := (g(t)^{-1}A, g(t)^{-1}u)$ satisfies (\ref{floweq2}) and its limit $(A_{\infty}, u_{\infty})$ solves $*F_{A_{\infty}} + \mu(u_{\infty})$ by Theorem \ref{ThmLimitStability}. By Proposition \ref{PropEngergy} and Theorem \ref{ThmLoj1} there exist $\gamma \in [\frac{1}{2}, 1)$ and $C,T > 0$ such that for all $t > T$
		\begin{align*}
				||*F_{g(t)^{-1}A} + \mu(g(t)^{-1}u) ||_{L^2}^2 &= 2\left( \mathcal{F}(A(t),u(t)) - \mathcal{F}(A_{\infty},u_{\infty}) \right) \\
																											 &= 2 \left( \mathcal{YMH}(A(t),u(t)) - \mathcal{YMH}(A_{\infty}, u_{\infty}) \right) \\
																											 &\leq 2 \left( \mathcal{YMH}(A(t),u(t)) - \mathcal{YMH}(A_{\infty}, u_{\infty}) \right)^{\gamma} \\
																											 &\leq C ||\nabla \mathcal{YMH}(A(t),u(t))||_{L^2} \\
																											 &= C ||\partial_t (A(t),u(t))||_{L^2}.
		\end{align*}
Theorem \ref{ThmConv} shows that the right-hand-side is integrable and hence
					$$m := \lim_{t\rightarrow \infty} \Psi_{(A,u)}(g(t)) > - \infty.$$ 
We claim $m = \inf \Psi_{(A,u)}$. For this let $\tilde{g}_0 \in (\mathcal{G}^c)^{2,2}(P)$ and denote by $\tilde{g}(t)$ the solution of (\ref{KNeq3}) starting at $\tilde{g}_0$. It follows from Step 1 of the proof of Theorem \ref{ThmNU}, that the pointwise geodesic distance between $g(t)$ and $\tilde{g}(t)$ in $G^c/G$ remains uniformly bounded. Since $\Psi_{(A,u)}$ is convex along geodesics in $(\mathcal{G}^c)^{2,2}(P)/\mathcal{G}^{2,2}(P)$ by (\ref{KNeq4}) and its gradient converges to zero along $g(t)$ and $\tilde{g}(t)$, it follows that $|\Psi_{(A,u)}(g(t)) - \Psi_{(A,u)}(\tilde{g}(t))|$ converges to zero. This proves the claim and $\Psi_{(A,u)}$ is bounded below $m$.
\end{proof}

\begin{Remark}[\textbf{The stable case}] \label{RmkStableKN}
In finite dimensions the Kempf--Ness functional of a point is proper if and only if this point is stable. Mundet \cite{Mundet:2000} established the following analogous result for the vortex in equations in great generality: $(A,u)$ is stable if and only if the complexified orbit $\mathcal{G}^c(A,u)$ has discrete $\mathcal{G}^c$-isotropy and for every $R > 0$ there exist $c_1, c_2 > 0$ such that
			\begin{align} \label{Meq1} ||*F_{e^{-\textbf{i}\xi}A} + \mu (e^{-\textbf{i}\xi} A)||_{L^2} < R \quad \Longrightarrow \quad 
				||\xi||_{L^{\infty}} \leq c_1\Psi_{(A,u)}(e^{\textbf{i}\xi}) + c_2. \end{align}
\end{Remark}

\section{Polystability and the moment-weight inequality}

\subsection{The Kobayashi--Hithchin correspondence}

\subsubsection*{Finite weights}

The weights of $(A,u) \in \mathcal{H}(P,X)$ are defined as the asymptotic slopes of $\Psi_{(A,u)}$ along the geodesic rays $[\exp(-\textbf{i}t\xi)]$ in $\mathcal{G}^c/\mathcal{G}$. Here $\xi \in \text{Lie}(\mathcal{G})$ is a section of $T^*\Sigma\otimes \text{ad}(P)$ and one may hope to replace the conditions on $\Psi_{(A,u)}$ in Theorem \ref{ThmKN} by conditions on these weights. In general, one needs to consider sections $\xi$ of very low regularity, namely of Sobolev class $H^1$. For bundles over a Riemann surface and smooth pairs $(A,u)$ every finite weight is obtained from a smooth section by Proposition \ref{PropFW} below.

\begin{Definition}
For $(A,u) \in \mathcal{H}(P,X)$ and $\xi \in H^1(\Sigma, \text{ad}(P))$ define
	$$w((A,u),\xi) := \lim_{t \rightarrow \infty} \langle *F_{e^{\textbf{i}t\xi} A} + \mu(e^{\textbf{i}t\xi} u), \xi \rangle_{L^2} \in \mathbb{R}\cup\{+\infty\}.$$	
By (\ref{KNeq4}), the right-hand-side is monotone increasing in $t$ and the limit exists.
\end{Definition}

Similarly, define by
		$$w(A,\xi) := \lim_{t\rightarrow \infty} \langle *F_{e^{\textbf{i}t\xi}A} , \xi \rangle, \qquad w(u,\xi) := \lim_{t\rightarrow \infty} \langle \mu(e^{\textbf{i}t\xi}u) , \xi \rangle$$
the weights for the $\mathcal{G}(P)$-action on $\mathcal{A}(P)$ and $\mathcal{S}(P,X)$ respectively. They are well-defined in $\mathbb{R}\cup\{+\infty\}$ and satisfy $w((A,u),\xi) = w(A, \xi) + w(u,\xi)$.

\begin{Proposition} \label{PropFW}
Let $A \in \mathcal{A}(P)$ be smooth and let $\xi \in H^1(\Sigma,\text{ad}(P))\backslash\{0\}$ with $w(A,\xi) < \infty$. 
		\begin{enumerate}
				\item Endow $P^c := P\times_G G^c$ with the holomorphic structure induced by $A$. Then there exists $\xi_0 \in \mathfrak{g}\backslash\{0\}$ and a holomorphic reduction $P_Q \subset P^c$ to the parabolic subgroup
							$$Q = Q(\xi_0) := \left\{ q \in G^c\,\left|\, \text{the limit $\lim_{t\rightarrow\infty} e^{\textbf{i}t\xi_0} q e^{-\textbf{i}t\xi_0} =: q_+$ exists} \right.\right\}.$$
				The reduction $P_Q \subset P^c$ induces a smooth reduction $P_K \subset P$ to the centralizer $K = C_G(\xi_0)$ and $\xi$ is the image of $\xi_0$ under the following map
								$$Z(\text{Lie}(K)) \rightarrow \Omega^0(\Sigma, \text{ad}(P_K)) \rightarrow \Omega^0(\Sigma,\text{ad}(P))$$
				where the first arrow identifies central elements with constant sections and the second map is obtained from the inclusion $P_K \subset P$.
				
				\item The limit $A_+ := \lim_{t\rightarrow\infty} e^{\textbf{i}t\xi} A$ exists in $H^1$ and $A_+$ restricts to a smooth connection on $P_K$.
		\end{enumerate}		
\end{Proposition}

\begin{proof}
This is an intrinsic version of \cite{Mundet:2000} Lemma 4.2 and makes use of a deep reularity result of Uhlenbeck and Yau \cite{UYau:1986} on weakly holomorphic subbundles.
See \cite{ST1} Lemma 5.7 for more details on the deduction. The reduction $P_K \subset P$ is induced by the isomorphism $G^c/Q(\xi_0) \cong G/C_G(\xi_0)$.
\end{proof}

\subsubsection*{Stable Kobayashi--Hitchin correspondence}

The Kobayashi--Hitchin correspondence for stable orbits says that $(A,u) \in \mathcal{H}(P,X)$ is stable if and only if $w((A,u), \xi) > 0$ for all $\xi \in \Omega^0(\Sigma, \text{ad}(P))\backslash\{0\}$. This was established by Mundet \cite{Mundet:2000} in greater generality and we briefly recall his argument. Suppose $(A,u)$ is stable and satisfies the vortex equation. Then
				\begin{align} \label{Seq1}
				w((A,u),\xi) = \langle *F_{A} + \mu(u), \xi \rangle_{L^2} + \int_0^{\infty} ||\mathcal{L}_{(e^{\textbf{i}t\xi}A, e^{\textbf{i}t\xi}u)} \xi||_{L^2}^2\, dt
		\end{align} 
is positive. It is a less obvious fact that this condition is $\mathcal{G}^c(P)$-invariant and hence $w(g(A,u),\xi) > 0$ for every $g \in \mathcal{G}^c(P)$. The converse direction depends on the Kempf--Ness theorem. Mundet shows by contradiction when no estimate (\ref{Meq1}) holds, then there exists a destabilizing direction $\xi$ with $w((A,u),\xi) \leq 0$. Once the estimate (\ref{Meq1}) is established, one obtains a solution to the vortex equation by direct methods of the calculus of variations. 

Our proof of the polystable case in Theorem \ref{ThmPS} below yields an alternative proof of the stable case under more restrictive assumptions.

\subsubsection*{Semistable Kobayashi--Hitchin correspondence}

We need to assume the following technical property for a pair $(A,u) \in \mathcal{H}(P,X)$:
		\begin{enumerate}[leftmargin=12mm]
				\item[\textbf{(H)}] If $\xi \in \Omega^0(\Sigma, \text{ad}(P))$ satisfies $w((A,u),\xi) \leq 0$ then $\sup_{t>0} ||\mu(e^{\textbf{i}t\xi} u)||_{L^2} < \infty.$
		\end{enumerate}	
We refer to Remark \ref{Rmk.H} for a discussion of this assumption. Following the ideas of Chen \cite{Chen:2009, Chen:2008}, Chen--Sun \cite{ChenSun:2010} and Donaldson \cite{Donaldson:2005} we prove the following version of the moment weight inequality which is Theorem \ref{ThmE} in the introduction.

\begin{Theorem}[\textbf{Sharp moment weight inequality}] \label{ThmSMWI}
Suppose $(A,u) \in \mathcal{H}(P,X)$ satisfies \textbf{(H)}. Then for all $\xi \in \Omega^0(\Sigma, \text{ad}(P))\backslash\{0\}$ it holds
		\begin{align} \label{MWIeq}  -\frac{w((A,u), \xi)}{||\xi||_{L^2}} \leq \inf_{g \in \mathcal{G}^c(P)} ||*F_{gA} + \mu(gu)||_{L^2}. \end{align}
If in addition \textbf{(A)}, \textbf{(B)}, \textbf{(C)} are satisfied and the right hand side is positive, then there exists a unique $\xi_0 \in \Omega^0(\Sigma, \text{ad}(P))$ with $||\xi||_{L^2} = 1$ which yields equality.
\end{Theorem}

\begin{proof} The proof is given in the next subsection on page \pageref{proofSMWI}.\end{proof}

\begin{Theorem}[\textbf{Semistable correspondence}] \label{ThmSS}
Assume \textbf{(A)}, \textbf{(B)}, \textbf{(C)} and suppose that $(A,u) \in \mathcal{H}(P,X)$ satisfies \textbf{(H)}. Then the following are equivalent:
		\begin{enumerate}
				\item $(A,u)$ is semistable in the sense of Definition \ref{DefStability}.
				\item $\inf_{g\in\mathcal{G}^c(P)} ||*F_{gA} + \mu(gu)||_{L^2} = 0$.
				\item $w((A,u),\xi) \geq 0$ for all $\xi \in \Omega^0(\Sigma, \text{ad}(P))$.
		\end{enumerate}
\end{Theorem}

\begin{proof}
This is a direct consequence of Theorem \ref{ThmNU} and Theorem \ref{ThmSMWI}.
\end{proof}

\subsubsection*{Polystable Kobayashi--Hitchin correspondence}

Consider for $(A,u) \in \mathcal{H}(P,X)$ the following properties
		\begin{enumerate}[leftmargin=12mm]
				\item[\textbf{(SS)}] For all $\xi \in \Omega^0(\Sigma,\text{ad}(P))$ it holds $w((A,u),\xi) \geq 0$.
		
				\item[\textbf{(PS1)}] For all $\xi \in \Omega^0(\Sigma,\text{ad}(P))$ with $\exp(\xi) = \mathds{1}$ and $(w(A,u),\xi) = 0$ the limit 
								$$\lim_{t\rightarrow \infty} e^{\textbf{i}t\xi}(A,u) \in (\mathcal{G}^c)^{2,2}(A,u)$$ 
							exists in $H^1\times H^2$ and remains in the (Sobolev completion of the) complexified group orbit $(\mathcal{G}^c)^{2,2}(A,u)$.
				\item[\textbf{(PS2)}] For all $\xi \in \Omega^0(\Sigma,\text{ad}(P))$ with $(w(A,u),\xi) = 0$ the limit 
								$$\lim_{t\rightarrow \infty} e^{\textbf{i}t\xi}(A,u) \in (\mathcal{G}^c)^{2,2}(A,u)$$ 
							exists in $H^1\times H^2$ and remains in the (Sobolev completion of the) complexified group orbit $\mathcal{G}^c(A,u)$.			
		\end{enumerate}

\begin{Theorem}[\textbf{Polystable correspondence}] \label{ThmPS}
Assume (\textbf{A}), (\textbf{B}), (\textbf{C}) and \textbf{(H)}. Then the following are equivalent
		\begin{enumerate}
				\item $(A,u)$ is polystable, i.e. there exits $g \in \mathcal{G}^c(P)$ such that $*F_{gA} + \mu(gu) = 0$.
				\item $(A,u)$ satisfies \textbf{(SS)} and \textbf{(PS1)}.
				\item $(A,u)$ satisfies \textbf{(SS)} and \textbf{(PS2)}.
		\end{enumerate}		
\end{Theorem}

\begin{proof}
See page \pageref{proofPS}.
\end{proof}

Assumption \textbf{(H)} is only needed for the application of Theorem \ref{ThmSS}. For twisted Higgs-bundles over Riemann surface a polystable Kobayashi-Hithchin correspondence was established by Garc\'{\i}a-Prada, Gothen and Mundet \cite{Garcia:2009} by different methods. We present a more general proof following the ideas of Chen--Sun \cite{ChenSun:2010}.

\subsection{Proof of the moment-weight inequality}
\label{proofSMWI} The purpose of this section is to prove Theorem \ref{ThmSMWI}. Section \ref{ProofSS1} contains the proof of the inequality (\ref{MWeq1}). The proof is essentially due to Chen \cite{Chen:2009, Chen:2008} and Donaldson \cite{Donaldson:2005}. Section \ref{ProofSS2} contains a proof of the equality in the unstable case. This is the analog of the Kempf existence theorem in finite dimension. The proof is based on arguments given by Chen--Sun \cite{ChenSun:2010} in the finite dimensional differentiable case. Section \ref{ProofSS3} contains a proof of the uniqueness claim. This is the analogue of the Kempf uniqueness theorem. The proof is the one given in \cite{RobSaGeo}, Theorem 11.3, for the finite dimensional setting and extends almost ad verbum to our setting.

\subsubsection{Proof of the inequality} \label{ProofSS1}
Let $(A,u) \in \mathcal{H}(P,X)$, $g_0 \in \mathcal{G}^c(P)$ and $\xi \in \Omega^0(\Sigma, \text{ad}(P))\backslash\{0\}$ be given and assume $w((A,u),\xi) \leq 0$. Define $\eta(t) \in \Omega^1(\Sigma, \text{ad}(P))$ and $u(t) \in \mathcal{G}(P)$ by
			\begin{align}  \label{MWeq1}g_0^{-1} = e^{-\textbf{i}\xi t} e^{-\textbf{i}\eta(t)} u(t). \end{align}
Let $\pi: G^c \rightarrow G^c/G$ denote the canonical projection. Since the left-invariant metric on $G^c/G$ has nonpositive curvature, the exponential map is distance increasing and it holds pointwise
		\begin{align*} \left|\left| \xi t - \eta(t) \right|\right| \leq \text{dist}_{G^c/G}(\pi(e^{\textbf{i}\xi t}), \pi(e^{\textbf{i}\eta(t)})) 
																								 \leq \text{dist}_{G^c/G}(\pi(\mathds{1}), \pi(g_0^{-1})). \end{align*}
In particular, there exists $C > 0$ such that $||\xi t - \eta(t) ||_{L^2} \leq C$ and this implies
		\begin{align} \label{MWeq2} \left|\left| \frac{\xi}{||\xi||_{L^2}} - \frac{\eta(t)}{||\eta(t)||_{L^2}}\right|\right|_{L^2} \leq \frac{C}{2 t ||\xi||_{L^2}} \end{align}
Define 
				$$g: [0,1] \rightarrow \mathcal{G}^c(P), \qquad g(s) := g_0^{-1} \exp \left(\textbf{i} s u(t)\eta(t) u^{-1}(t) \right).$$ 
Then $\gamma := \pi\circ g$ is the unique geodesic connecting $\pi(g_0^{-1})$ to $\pi(e^{-\textbf{i}t\xi})$. It follows from (\ref{KNeq1}) and the fact that $\Psi_{(A,u)}$ is convex along geodesics (\ref{KNeq4}) that
		\begin{align*}
				- ||*F_{g_0A} + \mu(g_0u)||_{L^2} &\leq \frac{1}{||\eta||_{L^2}} \alpha_{(A,u)}\left(\gamma(0), \dot{\gamma}(0) \right) \\
																					&\leq \frac{1}{||\eta||_{L^2}} \alpha_{(A,u)}\left(\gamma(1), \dot{\gamma}(1) \right) \\
																					&= \left\langle *F_{e^{\textbf{i}\xi t}A} + \mu(e^{\textbf{i}\xi t}u) , \frac{\eta(t)}{||\eta(t)||_{L^2}} \right\rangle_{L^2}.
		\end{align*}																			
Assumption \textbf{(H)}, Proposition \ref{PropFW} and (\ref{MWeq2}) show that the right-hand side converges to $\frac{w((A,u), \xi)}{||\xi||_{L^2}}$ and this completes the proof.

\subsubsection{Existence of the dominant weight} \label{ProofSS2}
Suppose that $(A_0,u_0) \in \mathcal{H}(P,X)$ is unstable. We prove in this section that there exists $\xi \in \Omega^0(\Sigma, \text{ad}(P))$ such that 
		\begin{align} \label{domeq} -\frac{w((A_0,u_0),\xi)}{||\xi||_{L^2}}  = \inf_{g\in\mathcal{G}^c} ||*F_{gA_0} + \mu(gu_0)||_{L^2}. \end{align}
Let $(A,u):[0,\infty) \rightarrow \mathcal{H}(P,X)$ be the solution of (\ref{floweq2}) starting at $(A_0,u_0)$, let $g : [0,\infty) \rightarrow \mathcal{G}^c(P)$ be the solution of (\ref{KNflow}). Define $\xi(t) \in \Omega^0(\Sigma, \text{ad}(P))$ and $k(t) \in \mathcal{G}(P)$ by
			\begin{align} \label{domeq0} g(t) = e^{-\textbf{i}\xi(t)}k(t). \end{align}
The strategy of the proof is to show that the limit
			\begin{align} \label{domeq1} \lim_{t\rightarrow \infty} \frac{\xi(t)}{t} =: \xi_{\infty} \end{align}
exists in $W^{1,p}$ and satisfies (\ref{domeq}).\\

\textbf{Step 1:} \textit{The limit (\ref{domeq1}) exists in $L^2$.}\\

Denote by $\pi : G^c \rightarrow G^c/G$ the canonical projection and let $\gamma := \pi\circ g$. Since $g^{-1}\dot{g} =\textbf{i}(*F_A + \mu(u))$ takes values in $\textbf{i}\mathfrak{g}$, it holds $\nabla_t \dot{\gamma} = d\pi(g) \textbf{i}(\partial_t (g^{-1}\dot{g}))$ and Theorem \ref{ThmConv} yields the estimate
	\begin{align} \label{domeq2} \int_T^{\infty} ||\nabla_t \dot{\gamma}(t)||_{L^2} \, dt \leq C \int_T^\infty ||\partial_t (A,u)||_{H^1}\, dt \leq C T^{-\epsilon}. \end{align}
Define $\gamma_t : [0,t] \rightarrow \mathcal{G}^c(P)/\mathcal{G}(P)$ by $\gamma_t(s) := \pi\left( e^{-\textbf{i}\frac{\xi(t)}{t} s} \right)$. Pointwise this is the geodesic segment connection $\pi(\mathds{1})$ to $\gamma(t)$ and we define 
	$$\rho_t(s) : \Sigma \rightarrow \mathbb{R}, \qquad  \rho_t(s) = \text{dist}_{G^c/G} (\gamma(s), \gamma_t(s)).$$
Since $G^c/G$ has nonpositve sectional curvature, there holds pointwise the estimate $\ddot{\rho}_t(s) \geq - ||\nabla \dot{\gamma}(s)||$ (see \cite{RobSaGeo} Appendix A). Hence (\ref{domeq2}) yields
	\begin{align} \label{domeq3} ||\dot{\rho}_t(s)||_{L^2} \leq \int_s^\infty ||\nabla_t \dot{\gamma}(t)||_{L^2} \leq C s^{-\epsilon} \end{align}
and integrating this estimate shows
	\begin{align}	\label{domeq4} ||\rho_t(s)||_{L^2} \leq \int_0^s ||\dot{\rho}_t(r)||_{L^2} \, dr \leq C s^{1-\epsilon}. \end{align}
Since the exponential map on $G^c/G$ is distance increasing, it follows pointwise for $0 < t_1 < t_2$
	\begin{align} \label{domeq5} \left|\left| \frac{\xi(t_1)}{t_1} - \frac{\xi(t_2)}{t_2} \right|\right| \leq  \left|\left|\dot{\gamma}_{t_1}(0) - \dot{\gamma}_{t_2}(0) \right|\right| \leq \frac{\rho_{t_2}(t_1)}{t_1}\end{align}
Now (\ref{domeq4}) and (\ref{domeq5}) show that $\frac{\xi(t)}{t}$ is a $L^2$-Cauchy sequence and the limit (\ref{domeq1}) exists in $L^2$.\\

\textbf{Step 2:} \textit{The limit (\ref{domeq1}) exists in $W^{1,p}$ for every $p \in (2,\infty)$.}\\

Let $\xi(t)$ be as in (\ref{domeq0}) and define
			$$R(t) := *F_{e^{\textbf{i}\xi(t)}A_0} - *F_{A_0} + \mu(e^{\textbf{i}\xi(t)} u_0) - \mu(u_0).$$
A similar calculation as in the proof of Theorem \ref{ThmNU} shows
		\begin{align*}
				2\left\langle R(t), \xi(t) \right\rangle 
						&=\Delta ||\xi(t)||^2 + 2 \int_0^1 \left( ||d_{e^{\textbf{i}s\xi(t)} A_0} \xi(t)||^2 + ||L_{e^{\textbf{i}s\xi(t)}u_0} \xi(t)||^2 \right)\, ds \\
						&\geq 2 ||\xi(t)|| \Delta ||\xi(t)||.
		\end{align*}
Thus $||\xi(t)|| : \Sigma \rightarrow [0,\infty)$ are positive functions satisfying $\Delta \xi(t) \leq ||R(t)||$ at points where $\xi(t) \neq 0$. An argument of Donaldson \cite{Donaldson:1987} (see \cite{Simpson:1988} Prop 2.1) using the mean-value property of harmonic functions shows that this implies an estimate
		\begin{align} \label{domeq01} ||\xi(t)||_{C^0} \leq C \left(1 + ||R(t)||_{L^2} + ||\xi(t)||_{L^1} \right). \end{align}
Since $(A(t),u(t))$ satisfies (\ref{floweq2}) and
				$$||*F_{e^{\textbf{i}\xi(t)}A_0} + \mu(e^{\textbf{i}\xi(t)}u_0)||_{L^2} = ||*F_{A(t)} + \mu(u(t))||_{L^2}$$
the term $||R(t)||_{L^2}$ is uniformly bounded and (\ref{domeq01}) simplifies to
			\begin{align} \label{domeq6} ||\xi(t)||_{C^0} \leq C (1 + ||\xi(t)||_{L^1}). \end{align}
In particular, $\frac{\xi(t)}{t}$ is uniformly bounded in $C^0$. Since $g(t)^{-1}\bar{\partial}_{A_0} g(t) = (A_0 - A(t))^{0,1}$ is uniformly bounded in $H^1$ and pointwise
					$$||g^{-1}(t) \bar{\partial}_{A_0} g(t)||^2 = ||e^{\textbf{i}\xi(t)} \bar{\partial}_{A_0} e^{-\textbf{i} \xi(t)}||^2 + ||(\bar{\partial}_{A_0} k(t)) k(t)^{-1}||^2$$
it follows that $e^{\textbf{i}\xi(t)}\bar{\partial}_{A_0} e^{-\textbf{i}\xi(t)}$ is uniformly bounded in $L^p$ for every $p \in (1,\infty)$. Now
					$$e^{\textbf{i}\xi(t)}\bar{\partial}_{A_0} e^{-\textbf{i}\xi(t)} = t e^{\textbf{i}\xi(t)/t}\bar{\partial}_{A_0} e^{-\textbf{i}\xi(t)/t}$$
implies that $e^{\textbf{i}\xi(t)/t}\bar{\partial}_{A_0} e^{-\textbf{i}\xi(t)/t}$ converges to zero in $L^p$ and by elliptic regularity the limit (\ref{domeq1}) exists in $W^{1,p}$. \\

\textbf{Step 3:} \textit{The limit $\xi_{\infty}$ defined by (\ref{domeq1}) yields equality in (\ref{MWIeq}).}\\

The Kempf--Ness functional (\ref{KNeq2}) satisfies $\Psi_{(A_0,u_0)}(\mathds{1}) = 0$, decreases along $\gamma(t)$ and is convex along geodesics. Hence $\Psi_{(A_0,u_0)}( e^{\textbf{i}s\xi(t)}) \leq 0$ for $0 < s < t$ and by continuity with respect to the $W^{1,p}$-topology, it takes nonpositive values along the geodesic ray $\gamma_{\infty}(t) := \pi\left( e^{\textbf{i}\xi_{\infty}t} \right)$. This implies $w((A_0,u_0), \xi_{\infty}) \leq 0$ and $\xi_{\infty}$ is smooth by Proposition \ref{PropFW}. Using again that $\Psi_{(A_0,u_0)}$ is convex along geodesics it follows
			\begin{align} \label{domeq8} \left| \Psi_{(A_0,u_0)}( \gamma(t)) - \Psi_{(A_0,u_0)}(\gamma_{\infty}(t)) \right| \leq M \cdot \text{dist}_{\mathcal{G}^c/\mathcal{G}} (\gamma(t), \gamma_{\infty}(t)) \end{align}
where $\text{dist}_{\mathcal{G}^c/\mathcal{G}}$ denotes the $L^2$-geodesic distance and 
			$$M := \sup_{t> 0} \max\left\{||*F_{g(t)^{-1}A_0} + \mu(g(t)^{-1}u_0)||_{L^2}, \,\,||*F_{e^{\textbf{i}\xi_{\infty} t}A_0} + \mu(e^{\textbf{i}t\xi}u_0)||_{L^2} \right\}.$$
which is finite by (\textbf{H}) and Proposition \ref{PropFW}. As $t \rightarrow \infty$ in (\ref{domeq4}) one obtains $\text{dist}_{\mathcal{G}^c/\mathcal{G}} (\gamma(t), \gamma_{\infty}(t)) \leq C t^{1-\epsilon}$ and hence
			\begin{align} \label{domeq9} \left| \Psi_{(A_0,u_0)}( \gamma(t)) - \Psi_{(A_0,u_0)}\gamma_{\infty}(t) \right| \leq C t^{1-\epsilon} .\end{align}
Then
		\begin{align*}
				-w((A_0,u_0),\xi_{\infty}) &= \lim_{t\rightarrow \infty} \frac{1}{t} \int_0^{t} \langle e^{-\textbf{i}\xi_{\infty} s} (A_0,u_0), \xi_{\infty} \rangle \, ds \\
																	&= \lim_{t\rightarrow \infty} \frac{\Psi_{(A_0,u_0)}(\gamma_{\infty}(t))}{t} \\
																	&= \lim_{t\rightarrow \infty} \frac{\Psi_{(A_0,u_0)}(\gamma(t))}{t} \\
																	&= \lim_{t\rightarrow \infty} \frac{1}{t}\int_0^t ||*F_{A(s)} + \mu(u(s))||_{L^2}^2\, ds \\
																	&= ||*F_{A_{\infty}} + \mu(u_\infty)||_{L^2}^2
		\end{align*}															
By Theorem \ref{ThmLimitStability} 
			$$||*F_{A_{\infty}} + \mu(u_\infty)||_{L^2} = \inf_{g \in \mathcal{G}^c} ||*F_{gA_0} + \mu(gu_0)|| =: m$$ 
and thus $- w((A_0,u_0),\xi_{\infty}) = m^2$. Now
		\begin{align*}
				||\xi_{\infty}||_{L^2} = \lim_{t\rightarrow \infty} \left|\left|\frac{\xi(t)}{t} \right|\right|_{L^2} 
															 \leq \lim_{t\rightarrow \infty} \frac{1}{t}\int_0^t ||\dot{\gamma}(s)||_{L^2}\, ds 
															 = ||*F_{A_{\infty}} + \mu(u_{\infty})||_{L^2} = m
		\end{align*}
shows $-w(A,\xi_{\infty})/||\xi_{\infty}||_{L^2} \geq m$ and the converse inequality follows from (\ref{MWIeq}).

\subsubsection{Uniqueness of the dominant weight} \label{ProofSS3}
Suppose that $(A,u) \in \mathcal{H}(P,X)$ is unstable and $\xi_0, \xi_1 \in \Omega^0(\Sigma, \text{ad}(P))$ satisfy $||\xi_1||_{L^2} = ||\xi_2||_{L^2} = 1$ and
				$$-w((A,u),\xi_1) = -w((A,u),\xi_2) = \inf_{g\in\mathcal{G}^c(P)} ||*F_{gA} + \mu(gu)||_{L^2} =: m > 0.$$
We prove in the following that this implies $\xi_1 = \xi_2$.\\

Define $\eta(t) \in \Omega^0(\Sigma, \text{ad}(P))$ and $k(t) \in \mathcal{G}(P)$ by
				\begin{align} \label{KUeq1} e^{-\textbf{i}t\xi_0} e^{\textbf{i}\eta(t)} = e^{-\textbf{i}t \xi_1} k(t).\end{align}
Let $\pi : G^c \rightarrow G^c/G$	denote the canonical projection and let $p(t) := \pi( e^{-\textbf{i}t\xi_0} e^{\textbf{i}\eta(t)/2})$ denote the midpoint between the geodesic rays spanned by $\xi_1$ and $\xi_2$. Since $G^c/G$ has nonpositive curvature, the exponential map (based at $p(t)$) is distance increasing and this yields
				\begin{align*} 
							d (1, p(t))^2 &\leq \frac{d(\pi(\mathds{1}), \pi(e^{\textbf{i}\xi_1 t}))^2 + d(\pi(\mathds{1}), \pi(e^{\textbf{i}\xi_2 t}))^2}{2} - \frac{d(\pi(e^{\textbf{i}\xi_1 t}), \pi(e^{\textbf{i}\xi_2 t}))^2}{4}\\
																					&\leq t^2 \left( 1 - \frac{||\xi_1 - \xi_2||_{L^2}}{4} \right)
				\end{align*}
where $d(\cdot, \cdot) = \text{dist}_{\mathcal{G}^c/\mathcal{G}}(\cdot, \cdot)$	denotes the $L^2$-geodesic distance. Since $\Psi_{(A,u)}$ is convex along geodesics, it follows $\Psi_{(A,u)}(p(t)) \leq - tm$ and hence
				\begin{align} \label{KUeq2} \frac{\Psi_{(A,u)}(p(t))}{d (1, p(t))} \leq \frac{-m}{\sqrt{1 - ||\xi_0 - \xi_1||_{L^2}/4}}. \end{align}
Denote for $r > 0$
			$$S_r := \left\{ \pi(e^{i\xi})\, \left| \,\xi \in \Omega^0(\Sigma, \text{ad}(P)), \, ||\xi||_{L^2} = r \right.\right\} \subset \mathcal{G}^c(P)/\mathcal{G}(P).$$
We claim
			\begin{align} \label{KUeq3} \lim_{r \rightarrow \infty} \frac{1}{r} \inf_{S_r} \Psi_{(A,u)} = -m .\end{align}
As $t \rightarrow \infty$ in (\ref{KUeq2}) the claim implies $\xi_1 = \xi_2$. The inequality "$\leq$" in (\ref{KUeq3}) follows by considering the values along the geodesics ray $\pi(e^{\textbf{i}\xi_1 t})$. For the other direction let $h \in \mathcal{G}^c(P)$ be given and using (\ref{KNeq4}) one estimates
					$$ \Psi_{h^{-1}(A,u)}(g) \geq - ||*F_{h^{-1}A} + \mu(h^{-1}u)||_{L^2} \cdot d(\pi(\mathds{1}), \pi(g))$$
Suppose $h$ is chosen such that  $||*F_{h^{-1}A} + \mu(h^{-1}u)||_{L^2} \leq m + \epsilon$. Then
			\begin{align*}
					\Psi_{(A,u)} (g) &= \Psi_{h^{-1}(A,u)}(h^{-1} g) + \Psi_{(A_0,u_0)}(h) \\
													 &\geq (-m-\epsilon) d(\pi(\mathds{1}), \pi(h^{-1}g)) + \Psi_{(A_0,u_0)}(h) \\
													 &\geq (-m-\epsilon) d(\pi(\mathds{1}), \pi(g)) + (-m-\epsilon) d(\pi(\mathds{1}), \pi(h^{-1})) + \Psi_{(A_0,u_0)}(h)
			\end{align*}
and as $d(\pi(\mathds{1}), \pi(g)) \rightarrow \infty$ and $\epsilon \rightarrow 0$ this proves (\ref{KUeq3}).

\subsection{Proof of the polystable correspondence}

The purpose of this section is to prove Theorem \ref{ThmPS}.

\begin{Proposition} \label{PropPS}
Let $(A,u) \in \mathcal{H}(P,X)$ be polystable, then $(A,u)$ satisfies \textbf{(SS)} and \textbf{(PS2)}.
\end{Proposition}

\begin{proof}[Proof of Proposition \ref{PropPS}]
Choose $g \in \mathcal{G}^c(P)$ such that $*F_{gA} + \mu(gu) = 0$. Then
		\begin{align} \label{PSeq0}
				w(g(A,u),\xi) = \langle *F_{gA} + \mu(gu), \xi \rangle_{L^2} + \int_0^{\infty} ||\mathcal{L}_{(e^{\textbf{i}t\xi}gA, e^{\textbf{i}t\xi}gu)} \xi||_{L^2}^2\, dt
		\end{align} 
shows $w(g(A,u),\xi) \geq 0$. Equality holds if and only if $\mathcal{L}_{g(A,u)} \xi = 0$ and $e^{\textbf{i}t\xi} g(A,u) = g(A,u)$ is constant. In particular, $g(A,u)$ satisfies \textbf{(SS)} and \textbf{(PS2)}. The Proposition follows now from Lemma \ref{LemmaPS} below.
\end{proof}

\begin{Lemma}\label{LemmaPS}
Let $(A,u) \in \mathcal{H}(P,X)$ and let $(B,v) \in \mathcal{G}^c(A,u)$.
		\begin{enumerate}
				\item If $(A,u)$ satisfies \textbf{(SS)} and \textbf{(PS1)} then $(B,v)$ satisfies \textbf{(SS)} and \textbf{(PS1)}.
				\item If $(A,u)$ satisfies \textbf{(SS)} and \textbf{(PS2)} then $(B,v)$ satisfies \textbf{(SS)} and \textbf{(PS2)}.
		\end{enumerate}
\end{Lemma}

\begin{proof}
We prove the second part first. Choose $g \in \mathcal{G}^c(P)$ such that $(B,v) = g(A,u)$ and let $\xi \in \Omega^0(\Sigma,\text{ad}(P))$ be such that $w(g(A,u), \xi) \leq 0$. Let $\xi_0 \in \mathfrak{g}$ and $P_Q \subset P^c$ be the $Q(\xi_0)$-bundle determined by $\xi$ as asserted in Proposition \ref{PropFW}. It is possible to decompose $g = q k$ with $q \in \mathcal{G}(P_Q)$ and $k \in\mathcal{G}(P)$ (e.g. by using the identity $G^c/B = G/Z(G)$ for any Borel subgroup $B \subset Q(\xi_0)$). By definition of $Q(\xi_0)$
				\begin{align} \label{PSeq1}
						q_+ := \lim_{t\rightarrow \infty} e^{\textbf{i}t\xi} q e^{-\textbf{i}t\xi}
				\end{align}		
Using the assumption $w(qk(A,u), \xi) \leq 0$ it follows for $t > 0$
	\begin{align} \label{PSeq2}
				0 \geq \Psi_{qk(A,u)}(e^{-\textbf{i}t\xi}) = \Psi_{k(A,u)} (q^{-1} e^{-\textbf{i}t\xi}) - \Psi_{k(A,u)}(q^{-1}) .
	\end{align}
Let $\pi: G^c \rightarrow G^c/G$ deonte the canoncial projection. Then
	\begin{align} \label{PSeq3}
			\text{dist}_{\mathcal{G}^c/\mathcal{G}}\left( \pi\left(e^{-\textbf{i}t\xi}\right), \pi\left(q^{-1} e^{-\textbf{i}t\xi}\right)\right) = \text{dist}_{\mathcal{G}^c/\mathcal{G}}\left(\pi(\mathds{1}), \pi\left(e^{\textbf{i}t\xi} q^{-1} e^{-\textbf{i}t\xi}\right)\right) \leq C
	\end{align}
which is bounded by (\ref{PSeq1}). For $t > s > 0$ define $\eta_{s,t} \in \Omega^0(\Sigma, \text{ad}(P))$ and $k_{s,t} \in \mathcal{G}(P)$ by
					$$e^{-\textbf{i}s\xi} e^{\textbf{i}\eta_{s,t}} = q^{-1} e^{-\textbf{i}t\xi} k_{s,t}.$$
Since the exponential map in $G^c/G$ is distance increasing, it follows from (\ref{PSeq3})
	\begin{align} \label{PSeq4}
			\lim_{t\rightarrow \infty} \left|\left| \frac{\eta_{s,t}}{t-s} - \xi\right|\right|_{L^2} = 0.
	\end{align}
If $\sup \{ \Psi_{k(A,u)}(e^{-\textbf{i}t\xi})\,|\, t > 0 \} < \infty$, then clearly $w(k(A,u), \xi) \leq 0$. Otherwise, (\ref{PSeq2}) shows that for all sufficently large $s > 0$ and every $t > s$ we have
		$$\Psi_{k(A,u)}(e^{-\textbf{i}s\xi}) > \Psi_{k(A,u)} (q^{-1} e^{-\textbf{i}t\xi}).$$ 
Since $\Psi_{k(A,u)}$ is convex along the geodesic segement $r \mapsto e^{-\textbf{i}s\xi} \cdot e^{\textbf{i}\eta_{s,t} r} $, it follows
					$$ d\Psi_{k(A,u)}\left( \pi(e^{-\textbf{i}s\xi}); d\pi(e^{-\textbf{i}s\xi})\textbf{i}\frac{\eta_{s,t}}{t-s} \right) = \left\langle *F_{e^{\textbf{i}s\xi}kA} + \mu (e^{-\textbf{i}s\xi}ku), \frac{\eta_{s,t}}{t-s} \right\rangle < 0. $$
Now (\ref{PSeq4}) implies $\left\langle F_{e^{\textbf{i}s\xi}A} + \mu(e^{\textbf{i}s\xi}u); \xi \right\rangle \leq 0$ for all sufficiently large $s$ and hence $w(k(A,u),\xi) \leq 0$.  Since $(A,u)$ satisfies \textbf{(SS)} by assumption, it follows $w((A,u),k^{-1}\xi k) = 0$ and $\textbf{(PS2)}$ implies that the limit
				$$(A_+,u_+) := \lim_{t\rightarrow \infty}  e^{\textbf{i}t k^{-1} \xi k } (A,u)$$
exists in $H^1\times H^2$ and $(A_+,u_+) \in (\mathcal{G}^c)^{2,2}(A,u)$. Hence
			$$(B_+,v_+) := \lim_{t\rightarrow \infty} e^{\textbf{i}t\xi} (B,v) = \lim_{t\rightarrow \infty} e^{\textbf{i}t\xi} q e^{-\textbf{i}t\xi} k e^{\textbf{i}t k^{-1} \xi k } (A,u) = q_+ k (A_+,u_+)$$
exists and $(B_+,v_+) \in (\mathcal{G}^c)^{2,2}(B,v)$. Moreover,
			$$w((B,v),\xi) = \langle*F_{q_+ k A_+} + \mu(q_+ k u_+),\xi \rangle$$ 
and it remains to verify that this vanishes. By Proposition \ref{PropFW}, there exists a reduction $P_K \subset P$ to the centralizer $K = C_G(\xi_0)$ and $q_+$ restricts to an element in $\mathcal{G}^c(P_K)$. Let $h: [0,1] \rightarrow \mathcal{G}^c(P_K)$ be a smooth path connecting $\mathds{1}$ to $q_+$ with $h^{-1}(t) \dot{h}(t) = \alpha(t) + \textbf{i}\beta(t)$ . A short calculation shows
		\begin{align*}
				& \partial_t \langle*F_{h(t) kA_+} + \mu(h(t) ku_+),\xi \rangle \\
				&\qquad =\left\langle - \left[ *F_{h(t)kA_+} + \mu(h(t)ku_+), \alpha(t) \right] , \xi \right \rangle
					+ \left\langle \mathcal{L}_{h(t)k(A_+,u_+)} \beta(t), \mathcal{L}_{h(t)k(A_+,u_+)} \xi \right\rangle \\
				&\qquad = 0
		\end{align*}		
where the last step uses $[\alpha(t),\xi] = 0$ and $\mathcal{L}_{h(t) k (A_+,u_+)} \xi = h(t)\mathcal{L}_{k(A,u)} \xi = 0$. Hence
			\begin{align*} w((B,v), \xi) &= \langle *F_{k A_+} + \mu(k u_+), \xi ) = \langle *F_{A_+}, \mu(u_+), k^{-1}\xi k \rangle \\
																	 &= w((A,u), k^{-1}\xi k) = 0
			\end{align*}														
and this completes the proof of the second part. 

The first part follows from the same argument, since $\exp(\xi) = \mathds{1}$ implies $\exp(k^{-1}\xi k) = \mathds{1}$.
\end{proof}

\begin{proof}[Proof of Theorem \ref{ThmPS}.] \label{proofPS}
If $(A,u)$ is polystable then it satisfies \textbf{(SS)} and \textbf{(PS2)} by Proposition \ref{PropPS}. For the converse direction let $(A_0,u_0) \in \mathcal{H}(P,X)$ be given and assume that it satisfies \textbf{(SS)} and \textbf{(PS1)}. Denote by $(A,u): [0,\infty) \rightarrow \mathcal{H}(P,X)$ the solution of (\ref{floweq2}) starting at $(A_0,u_0)$ with limit $(A_{\infty},u_{\infty})$. Theorem \ref{ThmSS} shows that $(A_0,u_0)$ is semistable and hence by Theorem \ref{ThmLimitStability} the limit solves $*F_{A_{\infty}} + \mu(u_{\infty}) = 0$. Denote the isotropy groups at the limit and their Lie algebras by
	$$H := \{ h \in \mathcal{G}^{2,2}(P) \,|\, (h A_{\infty}, h u_{\infty}) = (A_{\infty}, u_{\infty}) \}, \qquad \mathfrak{h} = \text{ker}(\mathcal{L}_{A_{\infty}, u_{\infty}})$$
	$$H^c := \{ h \in (\mathcal{G}^c)^{2,2}(P) \,|\, (h A_{\infty}, h u_{\infty}) = (A_{\infty}, u_{\infty}) \}, \qquad \mathfrak{h}^c = \text{ker}(\mathcal{L}^c_{(A_{\infty},u_{\infty})}).$$
By Theorem \ref{ThmLimitStability} we may assume that these groups are not discrete. There are two important properties to note: (1) By Lemma \ref{Lemma.SobolevReg} $A_{\infty}$ is gauge equivalent to a smooth connection. In particular, as a subgroup of the isotropy group of $A_{\infty}$, one can identify $H$ with a closed and hence compact subgroup of $G$. (2) Using the equation $*F_{A_{\infty}} + \mu(u_{\infty}) = 0$, a short calculation shows that $H^c$ is indeed the complexification of $H$.\\

\textbf{Step 1:} \textit{There exists an $H$-invariant holomorphic coordinate chart
				$$\psi: (T_{(A_{\infty},u_{\infty})} \mathcal{A}^{1,2}(P)\times \mathcal{S}^{2,2}(P,X), 0) \rightarrow \left(\mathcal{A}^{1,2}(P)\times \mathcal{S}^{2,2}(P,X),  (A_{\infty},u_{\infty}) \right)$$
defined on a neighborhood of the origin satisfying $d\psi(0,0) = \text{id}$.}\\

Let $\{g_p\}_{p\in P}$ be a smooth $G$-invariant family of Riemannian metrics on $X$, compatible with the holomorphic structure, such that $g_p$ is flat in a neighborhood of $u_{\infty}(p)$. Then $\exp_{g_p}: T_{u_{\infty}(p)} X \rightarrow X$ is holomorphic in a neighbourhood of the origin and $\phi_{\mathcal{S}}(\hat{u})(p) := \exp_{g_p}(\hat{u}(p))$ provides a holomorphic chart for $\mathcal{S}(P,X)$. Define
				$$\psi(a, \hat{u}) := \int_{H} h^{-1} (A_{\infty} +hah^{-1}) \, d\mu_H(h) 
								+ \phi_{\mathcal{S}} \int_H \phi_{\mathcal{S}}^{-1} \left( h^{-1} \phi_{\mathcal{S}}(h \hat{u})) \right)\,d\mu_H(h) $$
where $\mu_H$ denotes the Haar-measure on $H$ with $\mu_H(H) = 1$. This is well-defined for $||\hat{u}||_{L^{\infty}} \leq c ||\hat{u}||_{H^2}$ sufficiently small and satisfies the desired properties.\\

\textbf{Step 2:} \textit{The linearization of the holomorphicity condition $\bar{\partial}_A u = 0$ is the operator $D: H^{1}(\Sigma, \text{ad}(P)) \oplus H^2(\Sigma, u_{\infty}^*TX) \rightarrow H^1(\Sigma, \Lambda^{0,1} \otimes u_{\infty}^*TX)$
					$$D(a,\hat{u}) = \left(\nabla_{A_{\infty}} \hat{u} + L_{u_{\infty}} a\right)^{0,1}.$$
There exists an $H$-invariant holomorphic coordinate chart $\psi$ as in Step 1 with the additional property that 
				$$\psi(a,\hat{u}) \in \mathcal{H}^{1,2}(P,X) \qquad \Longrightarrow \qquad D(a,\hat{u}) = 0$$ 
for every pair $(a,\hat{u})$ in the domain of $\psi$.}\\
							
Since $\nabla_{A_{\infty}}$ is a Fredholm operator with closed range and finite dimensional cokernel, it follows that the image of $D$ is closed with finite codimension. Now any choice of complements for the kernel and image of $D$ yield a pseudoinverse 
				$$T: H^1(\Sigma, \Lambda^{0,1} \otimes u_{\infty}^*TX) \rightarrow H^{1}(\Sigma, \text{ad}(P)) \oplus H^2(\Sigma, u_{\infty}^*TX)$$
which is a bounded linear operator satisfying $DTD = D$ and $TDT = T$. Since $D$ is complex linear we can choose complex complements to obtain a complex linear pseudoinverse $T$. Moreover, $D$ is $H$-equivariant and for every $h \in H$ the operator $T_h := h T h^{-1}$ yields another complex linear pseudoinverse for $D$. The average
							$$\int_H \int_H T_{h_1} D T_{h_2} d\mu_H(h_1) d\mu_H(h_2)$$
with respect to the Haar measure $\mu_H$ provides a $H$-equivariant pseudoinverse.

Let $\psi$ be defined as in Step 1 and let $T$ be a complex linear $H$-invariant pseudoinverse of $D$. Consider on the domain of $\psi$ the map $\tilde{f}(a,\hat{u}) := f(\psi(a,\hat{u}))$
where $f(A,u) := \bar{\partial}_A u$. The map
				$$\theta(a,\hat{u}) := (a,\hat{u}) + T(\tilde{f}(a,\hat{u}) - D(a,\hat{u}))$$
satisfies $\theta(0) = 0$ and $d \theta(0) = \mathds{1}$. Hence there exists a local holomorphic inverse $\theta^{-1}$ around the origin by the inverse function theorem. It follows from the construction that $\tilde{f}_0 := (\mathds{1} - DT)\circ \tilde{f} \circ \theta^{-1}$ takes values in $\text{ker}(T)$ and
				$$\tilde{f}\circ\theta^{-1} = \tilde{f}_0 + D.$$
Since $\text{ker}(T)$ is a complement of $\text{Im}(D)$ this implies
					$$\tilde{f}\circ\theta^{-1}(a,\hat{u}) = 0 \qquad \Longleftrightarrow \qquad D(a,\hat{u}) = 0,\quad \tilde{f}_0(a,\hat{u}) = 0.$$
Step 2 follows from this discussion after replacing $\psi$ by $\psi\circ \theta^{-1}$.\\

\textbf{Step 3:} \textit{Denote by $\mathfrak{h}^\bot$ the $L^2$-orthogonal complement of $\mathfrak{h}$ in $H^2(\Sigma, \text{ad}(P))$ and by $V$ the $L^2$-orthogonal complement of the image of $\mathcal{L}^c_{A_{\infty},u_{\infty}}$. Then there exists $t_0 > 0$ and maps 
				$$(a,\hat{u}): [t_0, \infty) \rightarrow \text{ker}(D) \cap V, \qquad \xi,\eta: [t_0,\infty) \rightarrow \mathfrak{h}^\bot$$
such that $(a(t), \hat{u}(t))$ is in the domain of the chart $\psi$ constructed in Step 2 and
				\begin{align} \label{TPseq2} (A(t),u(t)) = e^{\textbf{i}\eta(t)} e^{\xi(t)} \psi(a(t), \hat{u}(t)) \end{align}
for all $t> t_0$.}\\

The map $\mathfrak{h}^\bot \times \mathfrak{h}^\bot\times V \rightarrow \mathcal{A}^{1,2}(P)\times \mathcal{S}^{2,2}(P,X)$ defined by 
				$(\xi, \eta, (a,\hat{u})) \mapsto e^{\textbf{i}\eta} e^{\xi} \psi(a,\hat{u})$
is smooth near the origin with invertible derivative. Step 3 follows now from the implicit function theorem and Step 2.\\

\textbf{Step 4:} \textit{Let $g: [0, \infty) \rightarrow \mathcal{G}^c(P)$ be the solution of the equation $g^{-1}\dot{g} = \textbf{i}(*F_{A(t)}+ \mu(u(t)))$ with $g(0) = \mathds{1}$ obtained in Theorem \ref{ThmVen}. There exists $t_1 \geq t_0$ with the following significance:
		$$h : [t_1, \infty) \rightarrow (\mathcal{G}^c)^{2,2}(P), \qquad h(t) := e^{-\textbf{i}\xi(t)} e^{-\textbf{i}\eta(t)} g^{-1}(t) g(t_1) e^{\textbf{i}\eta(t_1)} e^{\textbf{i}\xi(t_1)}$$	
satisfies $h(t) \in H^c$ and $(a(t),\hat{u}(t)) = h(t)^{-1}(a(t_1),\hat{u}(t_1))$ for every $t \geq t_1$.}\\

Let $t_1 \geq t_0$ be fixed. Rewrite the identity $(A(t),u(t)) = g(t)^{-1}(A_0,u_0)$ as $\psi(a(t), \hat{u}(t)) = h(t) \psi(a(t_1), \hat{u}(t_1))$ and differentiate this to obtain
		\begin{align} \label{TPSeq3} d \psi(a(t),\hat{u}(t))[\partial_t(a(t),\hat{u}(t))] = \partial_t (\psi(a(t), \hat{u}(t))) = \mathcal{L}^c_{\psi(a(t), \hat{u}(t))} \dot{h}(t)h^{-1}(t). \end{align}
For $(a,\hat{u}) \in V$ consider the operator
		$$ N_{(a,\hat{u})} : \, H^2(\Sigma, \text{ad}(P)^c) \times V \rightarrow T_{(A_{\infty}, u_{\infty})} ( \mathcal{A}^{1,2}(P)\times \mathcal{S}^{2,2}(P,X) )$$
				$$N_{(a, \hat{u})}(\zeta, (b, \hat{v})) := \mathcal{L}^c_{\psi(a,\hat{u})} \zeta - d \psi(a,\hat{u}) [ b, \hat{v} ].$$
Then (\ref{TPSeq3}) can be reformulated as
		\begin{align} \label{TPSeq4}  (h(t)^{-1}\dot{h}(t) , \partial_t (a(t),\hat{u}(t))) \in \text{ker}( N_{(a(t),\hat{u}(t))}) .\end{align}
Since $N_{(0,0)}$ is surjective with kernel $\mathfrak{h}^c$, it follows that $N_{(a,\hat{u})}$ is a surjective Fredholm operator with index $\text{dim}(\mathfrak{h}^c)$ for $||a||_{H^1} + ||\hat{u}||_{H^2}$ sufficiently small. For $\xi \in \mathfrak{h}$ it holds
					$$\left.\frac{d}{ds}\right|_{s=0} \psi(e^{\xi s} (a,\hat{u})) = \left.\frac{d}{ds}\right|_{s=0} e^{\xi s} \psi(a,\hat{u}) = \mathcal{L}_{\psi(a,\hat{u})} \xi.$$
Since $V$ is $H$-invariant this shows $\mathcal{L}_{\psi(a,\hat{u})} \mathfrak{h} \subset d\psi(a,\hat{u}) V$. Moreover, since $\psi$ is holomorphic and $V$ a complex subspaces, it follows $\mathcal{L}^c_{\psi(a,\hat{u})} \mathfrak{h}^c \subset d\psi(a,\hat{u}) V$. This implies that the kernel of $N_{(a,\hat{u})}$ projects onto $\mathfrak{h}^c$. For sufficiently large $t_1$ the same is true for all operators $N_{(a(t),\hat{u}(t))}$ with $t \geq t_1$. Then (\ref{TPSeq4}) shows $h^{-1}(t)\dot{h}(t) \in \mathfrak{h}^c$ for all $t \geq t_1$ and hence $h(t) \in H^c$. Since $\psi$ is holomorphic and $H$-equivariant this completes the proof of Step 4.
\\

\textbf{Step 5:} \textit{There exists $\xi_0 \in \mathfrak{h}$ with $\exp{\xi_0} = \mathds{1}$ such that $w(h(t_1)^{-1}(A(t_1),u(t_1)), \xi_0) = 0$ and
							\begin{align} \label{TPSeq5} \lim_{t\rightarrow \infty} e^{\textbf{i}t\xi_0} h(t_1)^{-1} (A(t_1), u(t_1)) = (A_{\infty}, u_{\infty}). \end{align}
In particular, $(A_{\infty},u_{\infty}) \in (\mathcal{G}^c)^{2,2}(A_0,u_0)$ and $(A_0,u_0)$ is polystable.}\\

The group $H$ acts on the finite dimensional vector space $X_0 := V \cap \text{ker}(D)$ by unitary automorphism. Step 4 shows that the origin is contained in the closure of the $H^c$-orbit of $(a(t_1),\hat{u}(t_1))$. The classical Hilbert-Mumford criterion (see \cite{RobSaGeo} Theorem 14.2) shows that there exists $\xi_0 \in \mathfrak{h}$ with $\exp{\xi_0} = \mathds{1}$ such that 
			$$\lim_{t\rightarrow \infty} e^{\textbf{i}t \xi_0} (a(t_1),\hat{u}(t_1)) = 0.$$ 
Since $\psi(e^{\textbf{i}t\xi_0}( a(t_1), \hat{u}(t_1)) = e^{\textbf{i}t\xi_0} \psi(a(t_1),\hat{u}(t_1)) = e^{\textbf{i}t\xi_0} h(t_1)^{-1} (A(t_1),u(t_1))$ for all $t \geq 0$, it follows
					$$\lim_{t\rightarrow \infty} e^{\textbf{i}t\xi_0} h(t_1)^{-1} (A(t_1), u(t_1)) = \psi(0) = (A_{\infty}, u_{\infty})$$
and $w(h(t_1)^{-1}(A(t_1),u(t_1)), \xi_0) = \langle *F_{A_{\infty}} + \mu(u_{\infty}), \xi_0 \rangle = 0$.

By Lemma \ref{Lemma.SobolevReg} there exists $k \in \mathcal{G}^{2,2}(P)$ such that $k(A_{\infty},u_{\infty})$ is smooth. Then $kh(t_1)^{-1} \in \text{ker}(\mathcal{L}^c_{k(A_{\infty},u_{\infty})}$ is smooth and 
				$$w(kh(t_1)^{-1}(A(t_1),u(t_1)), k\xi_0 k^{-1}) = w(h(t_1)^{-1}(A(t_1),u(t_1)), \xi_0) = 0$$
By Lemma \ref{LemmaPS} $kh(t_1)^{-1} (A(t_1),u(t_1))$ satisfies \textbf{(PS1)} and together with (\ref{TPSeq5}) this yields
	\begin{align*}	
		(A_{\infty},u_{\infty}) &= \lim_{t\rightarrow \infty} e^{\textbf{i}t\xi_0} h(t_1)^{-1} (A(t_1), u(t_1)) \\
														&= k^{-1} \lim_{t\rightarrow \infty} e^{k\xi_0 k^{-1}} k h(t_1)^{-1} (A(t_1),u(t_1)) \in (\mathcal{G}^c)^{2,2}(A_0,u_0)
	\end{align*}													
Hence $(A_0,u_0)$ is polystable by Theorem \ref{ThmLimitStability} and this completes the proof.

\end{proof}

\appendix
\section{The {\L}ojasiewicz inequality for Gelfand triples}

We provide an abstract version of the {\L}ojasiewicz gradient inequality following closely the arguments of R\r{a}de \cite{Rade1992} and Simon \cite{Simon:1983}.

Let $H$ be Hilbert space and let $V \subset H$ be a dense subset. Suppose $V$ is a Hilbert space in its own right with respect to an inner product $\langle\cdot, \cdot \rangle_V$ and assume that the inclusion $V \subset H$ is compact. Identifying $H$ with its dual, we obtain the Gelfand triple $V \subset H = H^* \subset V^*$. Let $F: V \rightarrow \mathbb{R}$ be a real analytic function and denote its differential by
			$$M := d F : V \rightarrow V^*.$$
Assume $F$ vanishes to the first order at the origin, i.e. $F(0) = 0$ and $M(0) = 0$. The linearization of $M$ at the origin is given by
			$$L = dM(0): V \rightarrow V^*$$
and we call this map the Hessian of $F$ at the origin.

\begin{Theorem} \label{ThmA}
Assume the setting described above and suppose there are constants $\delta, c > 0$ such that 
				\begin{align} \label{Aeq1} ||L x||_{V^*} \geq \delta ||x||_V - c ||x||_{H} \end{align}
is satisfied for all $x \in V$. Then there exist $\epsilon, C > 0$ and $\gamma \in [\frac{1}{2}, 1)$ such that for all $x \in V$ with $||x||_V \leq \epsilon$ it holds
					$$||dF(x)||_{V^*} \geq C |F(x)|^{\gamma}.$$
\end{Theorem}

\begin{proof}
The proof consists of six steps.\\

\textbf{Step 1:} \textit{$L$ has finite dimensional kernel and closed range.}\\

The proof is left as an exercise and uses the assumption that $V \subset H$ is compact. The result follows as in \cite{McSa:Jhol} Lemma A.1.1.\\

\textbf{Step 2:} \textit{Construction of the finite dimensional approximation.}\\

Let $K := \text{ker}(L)$ and denote its orthogonal complement by $W'$. The image $W'' := \text{Im}(L) \subset V^*$ agrees with the annihilator of $K$. Identifying $K^* \subset V^*$ with the annihilator of $W'$ yields decompositions 
				$$V = K \oplus W', \qquad V^* = K^*\oplus W'' $$
and $L$ restricts to an isomorphism $L: W' \rightarrow W''$. It follows from the implicit function theorem that there exists $\epsilon > 0$ and $\delta > 0$ such that for every $x \in K$ with $||x||_V < \epsilon$ there exists a unique $\phi(x) \in W'$ with $||\phi(x)||_V < \delta$ solving the equation $M(x + \phi(x)) \in K^*$. Moreover, the function
				$$\phi: B_{\epsilon}(0; K) \rightarrow B_{\delta}(0;W')$$
is analytic. Define
				$$f: B_{\epsilon}(0; K) \rightarrow \mathbb{R}, \qquad f(x) := F(x + \phi(x)).$$
This is a real analytic function on a finite dimensional domain.\\

\textbf{Step 3:} \textit{For $x \in B_{\epsilon}(0;K)$ it holds $d f (x) = M(x + \phi(x)) \in K^*$.}\\

For $x,y \in K$ the chain rule yields
			$$\langle df (x), y \rangle_{V^* \times V } = \langle M(x + \phi(x)), y + d \phi(x) y \rangle_{V^*\times V}.$$
and this proves the claim, since $M(x + \phi(x)) \in K^*$ annihilates $d \phi(x)y \in W'$.\\

\textbf{Step 4:} \textit{Decompose $x \in V$ with $||x||_V < \epsilon$ as 
					\begin{align} \label{Aeq2} x = x_0 + \phi(x_0) + x' \end{align}
with $x_0 \in K$ and $x' \in W'$. For sufficiently small $\epsilon > 0$ there exists $C > 0$ such that
				\begin{align} \label{Aeq3}  ||M(x)||_{V^*} \geq C\left(||df(x_0)||_{V^*} + ||x'||_V\right) 	\end{align}
holds for all $x \in V$ with $||x||_V < \epsilon$.				} \\

The terms in the decomposition (\ref{Aeq2}) satisfy the estimates
				\begin{align} \label{Aeq4} ||x_0||_V \leq C ||x||_V, \quad ||\phi(x_0)||_V \leq C ||x||_V, \quad ||x'||_V \leq C ||x||_V. \end{align}
Using Step 3 we obtain
		\begin{align*}
					M(x) &= M(x_0 + \phi(x_0) + x') \\
							 &= df(x_0) + \int_0^1 d M(x_0 + \phi(x_0) + tx')x' \, dt \\
							 &= df(x_0) + Lx' + \int_0^1 \left( d M(x_0 + \phi(x_0) + tx') - dM(0)\right) x'\, dt.
		\end{align*}
Since $dM$ is continuously differentiable, it follows from (\ref{Aeq4}) that there exists an estimate
				$$\sup_{t \in [0,1]} ||d M(x_0 + \phi(x_0) + tx') - dM(0)||_{\text{Hom}(V,V^*)} \leq C ||x||_V \leq C \epsilon.$$
Since $df(x_0) \in K^*$ and $L x' \in W''$ we have 
		$$||df(x_0) + L x'||_{V^*} \geq C \left(||df(x_0)||_{V^*} + ||L x'||_{V^*} \right) \geq C (||df(x_0)||_{V^*} + ||x'||_V).$$ 
Combining these estimates yields 
			$$||M(x)||_{V^*} \geq C_1(||df(x_0)||_{V^*} + ||x'||_{V}) - C_2 \epsilon ||x'||_{V}$$
and this proves (\ref{Aeq3}) after possibly shrinking $\epsilon > 0$. \\

\textbf{Step 5:} \textit{For sufficiently small $\epsilon > 0$ there exists $C > 0$ such that
															\begin{align} \label{Aeq5} |F(x)| \leq f(x_0) + C||x'||_V^2 \end{align} 
												for all $x \in V$ with $||x||_V < \epsilon$.}\\

The Taylor expansion of $F$ yields:
			\begin{align*}
					F(x) &= F(x_0 + \phi(x_0) + x')  \\
							 &= f(x_0) + \int_0^1 \langle M(x_0 + \phi(x_0) + t x'),  x' \rangle_{V^*\times V} \, ds \\
							 &= f(x_0) + \langle M(x_0 + \phi(x_0)), x' \rangle_{V^*\times V} \\
									& \qquad +\int_0^1 \int_0^1 \langle dM(x_0 + \phi(x_0) + st x') sx', x' \rangle_{V^*\times V}\, ds dt \\
							&=	f(x_0) + \langle df(x_0), x' \rangle_{V^*\times V} + \frac{1}{2} \langle L x', x' \rangle_{V^*\times V} + \langle L_2 x', x' \rangle_{V^*\times V}  \\	
			\end{align*}
where
		$$L_2 x' := \int_0^1 \int_0^1 s\left(dM(x_0 + \phi(x_0) + st x') - dM(0)\right) x' \, ds dt.$$
As in Step 4 one shows that this term satisfies an estimate 
				$$ \langle L_2 x', x' \rangle_{V^*\times V} \leq C ||x||_V||x'||_V^2 \leq C \epsilon ||x'||^2.$$
The open mapping theorem yields the estimate $\langle Lx' , x' \rangle_{V^*\times V} \geq C ||x'||_V^2$. Combining these estimates yields
							$$|F(x)| \leq f(x_0) + C_1 ||x'||_V^2 - C_2 \epsilon ||x'||_V^2$$
and this proves (\ref{Aeq5}) for sufficiently small $\epsilon > 0$.\\

\textbf{Step 6:} \textit{For suffiently small $\epsilon > 0$, there exists $C > 0$ and $\gamma \in[\frac{1}{2},1)$ such that
					\begin{align} \label{Aeq6} ||M(x)||_{V^*} \geq C |F(x)|^{\gamma} \end{align}
												 for all $x \in V$ with $||x||_V < \epsilon$.}\\

The gradient inequality of {\L}ojasiewicz \cite{Loj:1963} shows that for sufficently small $\epsilon > 0$ there exists $C > 0$ and $\gamma \in [\frac{1}{2}, 1)$ such that
									$$||df(x)||_V \geq |f(x)|^{\gamma}$$
for all $x \in K$	with $||x||_V < \epsilon$. Since $K$ is finite dimensional, there exists a constant such that $C ||df(x_0)||_{V^*} \geq  ||df(x_0)||_{V}$. Now the estimates (\ref{Aeq3}) and (\ref{Aeq5}) show
		\begin{align*}
					||M(x)||_{V^*}	&\geq C \left(||df (x_0) ||_{V^*} +||x'||_V\right) \\
													&\geq C_1 \left||F(x)| - C_2 ||x'||_V^2\right|^{\gamma} + C_3 ||x'||_V.
		\end{align*}
We may assume that $|F(x)| < 1$ for all $x \in B_{\epsilon}(0,V)$ and then follows (\ref{Aeq6}) with $C:=\min\{ C_1 2^{-\gamma}, C_3 /\sqrt{2 C_2} \}$.
\end{proof}

\newpage

\bibliographystyle{plain}
\bibliography{references}

\end{document}